# SOME COMPACTNESS CRITERIA FOR WEAK SOLUTIONS OF TIME FRACTIONAL PDES

LEI LI* AND JIAN-GUO LIU[†]

**Abstract.** The Aubin–Lions lemma and its variants play crucial roles for the existence of weak solutions of nonlinear evolutionary PDEs. In this paper, we aim to develop some compactness criteria that are analogies of the Aubin–Lions lemma for the existence of weak solutions to time fractional PDEs. We first define the weak Caputo derivatives of order $\gamma \in (0,1)$ for functions valued in general Banach spaces, consistent with the traditional definition if the space is $\mathbb{R}^d$ and functions are absolutely continuous. Based on a Volterra type integral form, we establish some time regularity estimates of the functions provided that the weak Caputo derivatives are in certain spaces. The compactness criteria are then established using the time regularity estimates. The existence of weak solutions for a special case of time fractional compressible Navier–Stokes equations with constant density and time fractional Keller–Segel equations in $\mathbb{R}^2$ are then proved as model problems. This work provides a framework for studying weak solutions of nonlinear time fractional PDEs.

**1. Introduction.** Memory effects are ubiquitous in physics and engineering, e.g. particles in heat bath ([1, 2]), soft matter with viscoelasticity ([3, 4]) can possess memory effects. Evolutionary equations of convolution type (see for examples [5, 6]) can be used to model these memory effects. When the memory effects have power law kernels, we can use fractional calculus to describe them [7, 8, 9, 10, 11]. There are two types of fractional derivatives that are commonly used: the Riemann–Liouville derivatives and the Caputo derivatives (See [9]). Caputo's definition of fractional derivatives was first introduced in [12] to study the memory effect of energy dissipation for some anelastic materials, and soon became a useful modeling tool in engineering. The use of Caputo derivatives can be justified, for example, from the generalized Langevin equation model [13], and by certain limiting processes and probability [14, 15]. Compared with Riemann–Liouville derivatives [7], Caputo derivatives remove singularities at the origin and share many similarities with the ordinary derivative so that they are suitable for initial value problems.

There are various definitions of Caputo derivatives in literature and they are all generalizations of the traditional Caputo derivatives. More recent definitions include those in [9, 16, 17, 18, 19]. In [9], the definition relies on Riemann–Liouville derivatives and is valid for some functions that do not necessarily have first derivatives; [16] relies on an integration by parts form and the functions only need to be Hölder continuous; in [17], some functional analysis approaches are used to extend the traditional Caputo derivatives to certain Sobolev spaces; the definition in [19] is based on the modified Riemann–Liouville operators and recovers the group structure. The underlying group structure mentioned in [19] is convenient for us to define the Caputo derivatives in even weaker spaces. In this paper, we will generalize the definition in [19] to weak Caputo derivatives for functions valued in general Banach spaces so that we can propose compactness criteria and study time fractional PDEs.

There is a significant amount of literature studying time-fractional ODEs (us-ing various definition of fractional derivatives) [8, 9, 10, 19, 20] and the theory is well-developed. Fractional stochastic differential equations have been discussed in [21, 22, 23, 13]. The fractional SDEs in [21, 22] are driven by fractional noise without fractional derivatives while the fractional SDEs in [23, 13] involve fractional derivatives. In [13], the authors argue that for physical systems, the derivatives paired with fractional Brownian noise must be Caputo derivatives following 'fluctuation-dissipation theorem'. In other situations (e.g. the finance model in [23]), Caputo derivatives and fractional Brownian motions may not be paired together.

*Department of Mathematics, Duke University, Durham, NC 27708, USA. (leili@math.duke.edu).
[†]Departments of Physics and Mathematics, Duke University, Durham, NC 27708, USA. (jliu@phy.duke.edu).

For time fractional PDEs with Caputo derivatives, the study is limited. The time fractional diffusion equations (using various definitions for Caputo derivatives) have been studied by many authors [24, 25, 26, 27, 28, 16, 29] and the theory is relatively well established. Time fractional Hamilton-Jacobi equations and the notion of viscosity solutions have been discussed in [30, 31] for examples. In [32], the author studied weak solutions of some linear evolutionary integro-differential equations, which include linear time fractional differential equations as special cases. Since the equation is linear, the weak solution can be obtained by weak compactness. The general nonlinear fractional PDEs are more challenging. One important class of methods used to study solutions of traditional nonlinear PDEs is to find some *a priori* estimates of Lyapunov functions and then to apply some strong compactness criteria. Some examples of (strong) compactness criteria include the Arzela–Ascoli theorem [33, Chap. 16], Kolmogorov-Riesz theorem [34], Rellich theorem [35], and the Aubin–Lions lemma [36, 37]. The Aubin–Lions lemma and its variants are for the compactness of functions over space and time, and are very useful for existence of weak solutions to nonlinear evolutionary PDEs. In this work, we aim to find suitable compactness criteria for nonlinear time fractional PDEs that are analogies of Aubin–Lions lemma (see Theorem 4.1 and Theorem 4.2) and see how these criteria can be used to prove the existence of weak solutions (Theorem 5.2 and Theorem 5.10).

The rest of the paper is organized as follows. In Section 2, we generalize the definition of Caputo derivatives to functions valued in general Banach spaces, and define the weak Caputo derivatives. Some basic properties of weak Caputo derivatives are then explored. In Section 3, we study the time regularity of functions from its weak Caputo derivatives. In particular, we have a time shift estimate for functions with weak Caputo derivatives in $L^p$. In Section 4, we establish the strong compactness criteria. Lastly in Section 5, we study a special case of time fractional compressible Navier–Stokes equations and time fractional Keller–Segel equations in $\mathbb{R}^2$ as model problems to show how the compactness criteria are used for the existence of weak solutions.

**2. Caputo derivatives based on a convolution group.** Before the general discussion, let us clarify some notations we use throughout this paper.

Let $X$ and $B$ be Banach spaces and $A \subset X$ is a connected Borel set in $X$. The notation $C_c^\infty(A;B)$ represents all the smooth functions $f : A \to B$ with compact supports. If the codomain $B$ is clear from the context, we may simply use $C_c^\infty(A)$ for short. (Usually $X = \mathbb{R}$ and $A = [a,b]$, $(a,b]$, $[a,b)$ or $(a,b)$). Similar notations are adopted for $L^p$ spaces and Sobolev spaces $W^{k,p}$. For example, $L^p((0,T);B)$ means all the Lebesgue measurable functions $f : (0,T) \to B$ so that $\int_0^T \|f\|_B^p dt < \infty$. If $B$ is clear, we may simply use $L^p(0,T)$. $f \in L_{\text{loc}}^p(A;B)$ means for any compact subset $K \subset A$, $f \in L^p(K;B)$. If $A$ is compact, then $L_{\text{loc}}^p$ is the same as $L^p$.

As we have mentioned, there are many definitions of Caputo derivatives, all reduce to the classical one if the function is smooth enough. In this paper, we work on mapping into general Banach spaces, so we must extend the definition of Caputo derivatives. We find that the approach in [19] is convenient for this generalization. For this purpose, we introduce more concepts and notations as follows.

DEFINITION 2.1. *Let $B$ be a Banach space. For a locally integrable function $u \in L^1_{\text{loc}}((0,T); B)$, if there exists $u_0 \in B$ such that*

$$\lim_{t \to 0+} \frac{1}{t} \int_0^t \|u(s) - u_0\|_B \, ds = 0, \tag{2.1}$$

*we call $u_0$ the right limit of $u$ at $t = 0$, denoted as $u(0+) = u_0$. Similarly, we define $u(T-)$ to be the constant $u_T \in B$ such that*

$$\lim_{t \to T-} \frac{1}{T-t} \int_t^T \|u(s) - u_T\|_B \, ds = 0. \tag{2.2}$$

The 'right limit' and 'left limit' should be understood in a weak sense. This is in fact the notion of Lebesgue point. Indeed, there may not exist a Lebesgue measure zero set $\mathcal{N}$ such that $\lim_{t \to 0+, t \notin \mathcal{N}} \|u(t) - u_0\|_B = 0$. As an example, consider $u : [0,1) \to \mathbb{R}$ such that $u(t) = n$ for $t \in [\frac{1}{n}, \frac{1}{n} + \frac{1}{n^4}]$ and $n \geq 2$. Otherwise, $u(t) = 0$. Then $u(0+) = 0$ in the sense of Definition 2.1, but $\lim_{t \to 0+, t \notin \mathcal{N}} |u(t) - 0| = 0$ does not hold for any Lebesgue measure zero set $\mathcal{N}$.

As in [19], we use the following distributions $\{g_\beta\}$ as the convolution kernels for $\beta > -1$:

$$g_\beta(t) := \begin{cases} \frac{\theta(t)}{\Gamma(\beta)} t^{\beta-1}, & \beta > 0, \\ \delta(t), & \beta = 0, \\ \frac{1}{\Gamma(1+\beta)} D\left(\theta(t) t^\beta\right), & \beta \in (-1, 0). \end{cases} \tag{2.3}$$

Here $\theta$ is the standard Heaviside step function, $\Gamma(\cdot)$ is the gamma function, and $D$ means the distributional derivative with respect to $t$.

$g_\beta$ can also be defined for $\beta \leq -1$ (see [19]) so that these distributions form a convolution group $\{g_\beta : \beta \in \mathbb{R}\}$, and consequently we have

$$g_{\beta_1} * g_{\beta_2} = g_{\beta_1 + \beta_2}, \tag{2.4}$$

where the convolution between distributions with one-sided bounded supports can be defined [38, Chap. 1].

**2.1. Functions valued in $\mathbb{R}^d$.** The fractional derivatives of a function valued in $\mathbb{R}^d$ can be defined by componentwise. Hence, it suffices to consider the derivatives for a function $u : (0, T) \to \mathbb{R}$, where $T \in (0, \infty]$.

In [19], the following modified Riemann–Liouville operators have been introduced

$$J_\alpha u := g_\alpha * (\theta u). \tag{2.5}$$

If $u$ is a distribution, $\theta u$ is defined as the weak limit (if it exists) of $\chi_n u$ as $n \to \infty$ where $\chi_n \in C_c^\infty(-1/n, \infty)$ is a smooth function that is 1 on $[0, \infty)$. The operators $\{J_\alpha : \alpha \in \mathbb{R}\}$ have group properties when acting on $u$ if $\theta u = u$ [19]. It is clear that if $\alpha > 0$ and $u \in L^1_{\text{loc}}(0, T)$, $J_\alpha$ is the fractional integral operator

$$J_\alpha u(t) = \frac{1}{\Gamma(\alpha)} \int_0^t (t-s)^{\alpha-1} u(s) \, ds. \tag{2.6}$$

DEFINITION 2.2 ([19]). *Let $0 < \gamma < 1$ and $u \in L^1_{\text{loc}}[0,T)$. Given any $u_0 \in \mathbb{R}$, the $\gamma$th order generalized Caputo derivative of $u$ from $t = 0$ associated with $u_0$ is a distribution in $\mathscr{D}'(-\infty, T)$ with support in $[0, T)$, given by*

$$(2.7) \qquad D_c^\gamma u := J_{-\gamma} u - u_0 g_{1-\gamma} = g_{-\gamma} * \Big(\theta(t)(u - u_0)\Big).$$

*If $u(0+)$ exists in the sense of Definition 2.1 and $u_0 = u(0+)$, we call $D_c^\gamma u$ the Caputo derivative of $u$.*

**Remark 2.1.** The generalized Caputo derivatives depends on the choice of $u_0$. For example, $u(t) = 1$. If we choose $u_0 = 1$, the weak Caputo derivative is zero while if we choose $u_0 = 0$, the generalized Caputo derivative is $\frac{\theta(t)}{\Gamma(1-\gamma)} t^{-\gamma}$. $t^{-\gamma}$ is like the Dirac delta for first derivative (if $f(t) = 1 + t$ and we choose $f_0 = 0$, then the first derivative becomes $\delta(t) + 1$ while choosing $f_0 = 1$ yields that $f' = 1$).

From here on, without mentioning, we always use $D_c^\gamma u$ to mean the Caputo derivative if $u(0+)$ exists (i.e. $u_0 = u(0+)$).

**Remark 2.2.** If $T < \infty$, $g_{-\gamma} * u$ should be understood as the restriction of the convolution onto $\mathscr{D}'(-\infty, T)$. One can refer to [19] for the technical details.

**Remark 2.3.** If there is a version of $u$ that is absolutely continuous on $(0, T)$ (still denoted as $u$), then the Caputo derivative is reduced to the traditional definition of Caputo derivative [19]:

$$(2.8) \qquad D_c^\gamma u(t) = \frac{1}{\Gamma(1-\gamma)} \int_0^t \frac{\dot{u}(s)}{(t-s)^\gamma} ds,$$

*where $\dot{u}$ means the time derivative of $u$.*

Definition 2.2 is more useful than the traditional definition (Equation (2.8)) (see for instance [8, 9, 10, 25]) theoretically, since it reveals the underlying group structure. With the assumption that $u$ is locally integrable and has a right limit at $t = 0$, Definition 2.2 and the group property (2.4) reveal that

$$(2.9) \qquad u(t) = u_0 + J_\gamma(D_c^\gamma u)(t) = u_0 + \frac{1}{\Gamma(\gamma)} \int_0^t (t-s)^{\gamma-1} D_c^\gamma u(s)\, ds.$$

Note that the integral simply means convolution $g_\gamma * D_c^\gamma u$. If $D_c^\gamma u \in L^1_{\text{loc}}[0,T)$, it can be understood in the Lebesgue integral sense. This simply means that the fractional integral of Caputo derivative recovers the function, so it is a fractional version of fundamental theorem of calculus. Consequently, we conclude that

LEMMA 2.3. *Suppose $E(\cdot) : [0, \infty) \to \mathbb{R}$ is continuous. If there exists $f(t) \in L^1_{\text{loc}}([0, \infty), \mathbb{R})$ such that on $(0, \infty)$*

$$D_c^\gamma E \le f,$$

*where this inequality means that $f - D_c^\gamma E$ is a non-negative distribution on $(0, \infty)$ (see [19]), then*

$$(2.10) \qquad E(t) \le E(0) + \frac{1}{\Gamma(\gamma)} \int_0^t (t-s)^{\gamma-1} f(s)\, ds, \ a.e.$$

Another property regarding convex functional that is important to us is as follows. For a detailed discussion, one can refer to [19, Proposition 3.11]. For the convenience of the readers, we provide a simple version and its concise proof here.

LEMMA 2.4. *If $u : [0, T] \to \mathbb{R}^d$ is $C^1((0, T); \mathbb{R}^d) \cap C([0, T]; \mathbb{R}^d)$, and $u \mapsto E(u)$ is a $C^1$ convex function on $\mathbb{R}^d$, then the Caputo derivative is given by*

$$(2.11) \qquad D_c^\gamma u(t) = \frac{1}{\Gamma(1-\gamma)} \left( \frac{u(t) - u(0)}{t^\gamma} + \gamma \int_0^t \frac{u(t) - u(s)}{(t-s)^{\gamma+1}} ds \right)$$

*and*

$$(2.12) \qquad D_c^\gamma E(u(t)) \leq \nabla_u E(u(t)) \cdot D_c^\gamma u.$$

*Proof.* The first claim follows from integration by parts of (2.8). For the second one, we note

$$E(u(t)) - E(b) \leq \nabla_u E(u(t)) \cdot (u(t) - b), \quad \forall b \in \mathbb{R}^d,$$

since $E(\cdot)$ is a convex function. Combining with the fact that $E(u(\cdot)) \in C^1((0, T); \mathbb{R}) \cap C([0, T]; \mathbb{R})$ (thus absolutely continuous), we have

$$D_c^\gamma E(u(t)) = \frac{1}{\Gamma(1-\gamma)} \left( \frac{E(u(t)) - E(u(0))}{t^\gamma} + \gamma \int_0^t \frac{E(u(t)) - E(u(s))}{(t-s)^{\gamma+1}} ds \right)$$
$$\leq \nabla_u E(u(t)) \cdot D_c^\gamma u. \qquad \square$$

Now, we move onto the right derivatives and integration by parts for fractional derivatives. In [19], there is another group given by

$$\tilde{\mathscr{C}} := \{\tilde{g}_\alpha : \tilde{g}_\alpha(t) = g_\alpha(-t), \alpha \in \mathbb{R}\}.$$

Clearly, $\operatorname{supp} \tilde{g} \subset (-\infty, 0]$. For $\gamma \in (0, 1)$:

$$(2.13) \qquad \tilde{g}_{-\gamma}(t) = -\frac{1}{\Gamma(1-\gamma)} D(\theta(-t)(-t)^{-\gamma}) = -D\tilde{g}_{1-\gamma}(t),$$

where $D$ means the distributional derivative on $t$. Suppose $\phi$ is absolutely continuous and $\phi = 0$ for $t > T$, then it can be verified directly that

$$\tilde{g}_{-\gamma} * \phi(t) = -\frac{1}{\Gamma(1-\gamma)} \frac{d}{dt} \int_t^\infty (s-t)^{-\gamma} \phi(s)\, ds$$
$$= -\frac{1}{\Gamma(1-\gamma)} \frac{d}{dt} \int_t^T (s-t)^{-\gamma} \phi(s)\, ds.$$

By the definition of $\tilde{g}_\alpha$, we have

LEMMA 2.5. *Suppose $\phi_1, \phi_2$ are distributions in $\mathscr{D}'$, and there exist $t_1, t_2 \in \mathbb{R}$ such that $\operatorname{supp} \phi_1 \subset (-\infty, t_1)$ while $\operatorname{supp} \phi_2 \subset (t_2, \infty)$. Then, if one of them is in $C_c^\infty(\mathbb{R})$ or both are absolutely continuous so that both $\langle g_{-\gamma} * \phi_1, \phi_2 \rangle$ and $\langle \phi_1, \tilde{g}_{-\gamma} * \phi_2 \rangle$ are defined, then it holds that*

$$(2.14) \qquad \langle g_{-\gamma} * \phi_1, \phi_2 \rangle = \langle \phi_1, \tilde{g}_{-\gamma} * \phi_2 \rangle.$$

Using the group $\tilde{\mathscr{C}}$, we define the right Caputo derivative as

DEFINITION 2.6. *Let $0 < \gamma < 1$. Consider $u \in L^1_{\text{loc}}(-\infty, T]$. Given $u_T \in \mathbb{R}$, the $\gamma$th order generalized right Caputo derivative (up to $T$) of $u$ associated with $u_T$ is a distribution in $\mathscr{D}'(\mathbb{R})$ with support in $(-\infty, T]$, given by*

$$(2.15) \qquad \tilde{D}_{c;T}^\gamma u := \tilde{g}_{-\gamma} * (\theta(T - t)(u(t) - u_T)).$$

*If $u$ has a left limit at $t = T$ in the sense of Definition 2.1 and $u_T = u(T-)$, we call $\tilde{D}_{c;T}^\gamma u$ the right Caputo derivative (up to $T$).*

From here on, without mentioning, we always use $\tilde{D}^\gamma_{c;T} u$ to mean the right Caputo derivative with the natural terminal value $u(T-)$ if this limit exists.

**Remark 2.4.** *If $u \in L^1_{\mathrm{loc}}(a, T), a < T$ and $u$ has a left limit $u(T-)$ at $t = T$, we do similar trick as that for $D^\gamma_c u$ in [19]; i.e. we extend $u$ to $(-\infty, T]$ by considering $u\chi_n$ where the smooth function $\chi_n$ is supported in $(a + \frac{1}{2n}, \infty)$ and equals one on $(a + 1/n, \infty)$. The weak limit of $\tilde{D}^\gamma_{c;T}(u\chi_n)$ as $n \to \infty$ in $\mathscr{D}'(a, T)$ is defined to be the right Caputo derivative.*

Similar as in [19], we can show that

**LEMMA 2.7.** *If $u$ is absolutely continuous on $(a, T)$, $a < T$, then*

$$\text{(2.16)} \qquad \tilde{D}^\gamma_{c;T} u(t) = -\frac{1}{\Gamma(1-\gamma)} \int_t^T (s-t)^{-\gamma} \dot{u}(s)\, ds, \ \forall t \in (a, T).$$

*Consequently, if $\varphi \in C_c^\infty(-\infty, T)$,*

$$\tilde{D}^\gamma_{c;T} \varphi(t) = \frac{-1}{\Gamma(1-\gamma)} \int_t^T (s-t)^{-\gamma} \dot{\varphi}(s)\, ds$$
$$= \tilde{g}_{-\gamma} * \varphi(t)$$
$$= \frac{-1}{\Gamma(1-\gamma)} \frac{d}{dt} \int_t^T (s-t)^{-\gamma} \varphi(s)\, ds.$$

We only sketch the proof while leave the details for the readers here. For the first claim, we note that if $u$ is absolutely continuous, then $\dot{u} \in L^1(a, T)$. The distributional derivative of $\theta(T-t)(u(t) - u(T-))$ is simply $\theta(T-t)\dot{u}$. Using the fact $\tilde{g}_{-\gamma} = -D\tilde{g}_{1-\gamma}$ yields the claim. For the second claim, we only have to justify $\tilde{D}^\gamma_{c;T}\varphi(t) = \tilde{g}_{-\gamma} * \varphi(t)$, which follows from the fact $\theta(T-t)(\varphi - \varphi(T)) = \varphi$ if $\varphi \in C_c^\infty(-\infty, T)$.

Using Lemma 2.5 and Definition 2.6, it is easy to obtain the integration by parts formula (we omit the proof as well):

**LEMMA 2.8.** *Let $u, v$ be absolutely continuous on $(0, T)$, then we have the integration by parts formula for Caputo derivatives*

$$\text{(2.17)} \qquad \int_0^T (D^\gamma_c u(t))(v(t) - v(T-))\, dt = \int_0^T (u(t) - u(0+))(\tilde{D}^\gamma_{c;T} v(t))\, dt.$$

*This relation also holds if $u \in L^1_{\mathrm{loc}}(0, T)$, $u(0+)$ is replaced with the assigned initial value $u_0$, and $v \in C_c^\infty(-\infty, T)$.*

**Remark 2.5.** *If $\gamma \to 1$, it is not hard to see that $\tilde{D}^\gamma_{c;T} u \to -u'(t)$ weakly. Hence, the right derivatives carry a natural negative sign.*

This identity is indeed not new under the classical definition [8, 9]. Using the convolution groups, the identity here becomes very natural and holds for larger class of functions.

**Remark 2.6.** *For Lemma 2.8, it might be illustrative to write out the computation for smooth $u$ and $v$ using traditional definitions:*

$$\int_0^T (D^\gamma_c u)(v(t) - v(T))dt = \int_0^T \frac{1}{\Gamma(1-\gamma)} \int_0^t \frac{\dot{u}(s)}{(t-s)^\gamma} ds\, (v(t) - v(T))dt$$
$$= \int_0^T \frac{\dot{u}(s)}{\Gamma(1-\gamma)} \int_s^T \frac{v(t) - v(T)}{(t-s)^\gamma} dt\, ds = -\int_0^T \frac{u(s) - u(0)}{\Gamma(1-\gamma)} \frac{d}{ds} \int_s^T \frac{v(t) - v(T)}{(t-s)^\gamma} dt\, ds.$$

*Moreover,*

$$-\frac{d}{ds}\int_s^T \frac{v(t)-v(T)}{(t-s)^\gamma}dt = -\frac{d}{ds}\int_0^{T-s}\frac{v(t+s)-v(T)}{t^\gamma}dt = -\int_0^{T-s}\frac{\dot{v}(t+s)}{t^\gamma}dt.$$

*Hence, the identity is verified.*

**2.2. Functions valued in general Banach spaces.** Now, we move onto Caputo derivatives for functions valued in general Banach spaces. We first define the weak Caputo derivatives as abstract linear functionals from $C_c^\infty((-\infty,T);\mathbb{R})$ to a Banach space $B$ for theoretical purposes. These functionals can be understood as the generalization of distributions (if $B=\mathbb{R}$, they are reduced to the usual distributions, as studied already in [19]). However, for practical purposes, we care more about those functions such that the Caputo derivatives are in $L^1_{\text{loc}}([0,T);B)$, for which we identity the action of the functionals with the $L^2$ pairing (see Equation (2.25)).

We now fix $T > 0$ and introduce the following set:

(2.18) $\quad \mathscr{D}' := \Big\{ v \mid v : C_c^\infty((-\infty,T);\mathbb{R}) \to B \text{ is a bounded linear operator} \Big\}.$

In other words, $\mathscr{D}'$ consists of functionals from $C_c^\infty((-\infty,T);\mathbb{R})$ to $B$. This is the analogy of the distributions $\mathscr{D}'(\mathbb{R})$ used in [19]. For $v \in \mathscr{D}'$, we say $\operatorname{supp} v \subset K$ for $K$ that is a closed subset of $(-\infty,T)$ (under the topology of $(-\infty,T)$), if for any $\varphi \in C_c^\infty((-\infty,T);\mathbb{R})$ with $\operatorname{supp}\varphi \subset (-\infty,T) \setminus K$, we have

$$\langle v, \varphi \rangle = 0.$$

With this notion, we can introduce

(2.19) $\quad\quad\quad\quad \mathscr{D}'_+ := \{v \in \mathscr{D}' \mid \operatorname{supp} v \subset [0,T)\}.$

Motivated by the usual weak derivatives of the functions valued in Banach spaces ([35, sect. 5.9.2]) and the integration by parts formula (Lemma 2.8), we define

DEFINITION 2.9. *Let $B$ be a Banach space and $u \in L^1_{\text{loc}}([0,T);B)$. Let $u_0 \in B$. We define the weak Caputo derivative of $u$ associated with initial data $u_0$ to be $D_c^\gamma u \in \mathscr{D}'$ such that for any test function $\varphi \in C_c^\infty((-\infty,T);\mathbb{R})$,*

(2.20) $\quad \langle D_c^\gamma u, \varphi \rangle := \int_{-\infty}^T (u-u_0)\theta(t)(\tilde{D}_{c;T}^\gamma \varphi)\,dt = \int_0^T (u-u_0)\tilde{D}_{c;T}^\gamma \varphi\,dt.$

*If $u(0+) = u_0$ in the sense of Definition 2.1 under the norm of the underlying Banach space $B$, we call $D_c^\gamma u$ the Caputo derivative.*

**Remark 2.7.** *Similar as in Remark 2.1, the choice of $u_0$ affects $D_c^\gamma u$. If $u(0+)$ exists, choosing $u_0 = u(0+)$ kills the singularity brought by the jump at $t=0$, and we will make this convention default. We also remark that the idea of using pairing to define fractional derivatives also appeared in [18].*

We have the following observation.

LEMMA 2.10. $D_c^\gamma u \in \mathscr{D}'_+$. *In other words*, $\operatorname{supp} D_c^\gamma u \subset [0,T)$.

*Proof.* By the explicit formula $\tilde{D}_{c;T}^\gamma \varphi(t) = -\frac{1}{\Gamma(1-\gamma)}\int_t^T (s-t)^{-\gamma}\dot\varphi(s)ds$, we find that if $\operatorname{supp}\varphi \subset (-\infty,0)$, the integral in Definition 2.9 is zero. $\square$

We now check that Definition 2.9 is consistent with the usual definitions.

LEMMA 2.11. *If $B = \mathbb{R}^d$, then the weak Caputo derivative in Definition 2.9 agrees with Definition 2.2.*

*Proof.* We only have to focus on $d = 1$, because for general $d$ we define them componentwise. Take $\varphi \in C_c^\infty((-\infty, T); \mathbb{R})$. By Lemma 2.7,
$$\tilde{D}_{c;T}^\gamma \varphi = \tilde{g}_{-\gamma} * \varphi.$$

The claim then follows from
$$\langle D_c^\gamma u, \varphi \rangle = \int_0^T (u - u_0) \tilde{D}_{c;T}^\gamma \varphi \, dt$$
$$= \int_0^T (u - u_0)(\tilde{g}_{-\gamma} * \varphi) \, dt$$
$$= \langle g_{-\gamma} * ((u - u_0)\theta(t)), \varphi \rangle. \qquad \square$$

Now, we investigate the properties of weak Caputo derivatives. The proof in Lemma 2.11 actually motivates us to consider the convolution between $g_{-\gamma}$ and distributions in $\mathscr{D}'$. Let $v \in \mathscr{D}'$ with $\operatorname{supp} v \subset [0, T)$. Consider a sequence of smooth functions $\chi_n$ that is 1 on $(-n, T - \frac{1}{n})$ and zero on $[T - \frac{1}{2n}, +\infty)$. Then, $\chi_n v$ is a distribution for $\varphi \in C_c^\infty(\mathbb{R}; \mathbb{R})$. Then, $g_\alpha * (\chi_n v)$ can be defined as in [19, Def. 2.1] and becomes a functional in $\mathscr{D}'$.

DEFINITION 2.12. *We define the convolution between $v$ and $g_\alpha$ as $g_\alpha * v \in \mathscr{D}'$:*

(2.21) $$g_\alpha * v := \lim_{n \to \infty} g_\alpha * (\chi_n v) \text{ in } \mathscr{D}'.$$

Similar as the cases for $\mathbb{R}^d$ discusses in [19], we define the following

DEFINITION 2.13. $J_\alpha : \mathscr{D}'_+ \to \mathscr{D}'_+$ *as*

(2.22) $$J_\alpha v := g_\alpha * v.$$

*When $\alpha > 0$, we call $J_\alpha$ the fractional integral operator with order $\alpha$.*

**Remark 2.8.** *For a general distribution $w \in \mathscr{D}'$, $J_\alpha w$ is defined to be $J_\alpha(\theta w)$ where $\theta w$ is the restriction of $w$ onto $\mathscr{D}'_+$ if it exists (see [19, Def. 2.14]). We will not consider this general case in this paper.*

The first property is the fractional fundamental theorem of calculus for functions valued in general Banach spaces.

LEMMA 2.14. *Let $u \in L^1_{\mathrm{loc}}((0, T); B)$. Assume the weak Caputo derivative for an assigned initial value $u_0$ is $D_c^\gamma u$. We have in $\mathscr{D}'$ that*

(2.23) $$(u - u_0)\theta(t) = J_\gamma(D_c^\gamma u) = g_\gamma * D_c^\gamma u.$$

*Proof.* To be convenient, we extend $u$ to be defined on $[0, \infty)$ by defining its values to be zero outside $[0, T)$.

We now pick $\eta \in C_c^\infty(0, 1)$, $0 \leq \eta \leq 1$ and $\int_{(0,1)} \eta \, dt = 1$. We define $\eta_\epsilon = \frac{1}{\epsilon} \eta\left(\frac{t}{\epsilon}\right)$. We note first that $\eta_\epsilon * [(u - u_0)\theta(t)]$ is zero at $t = 0$ and supported in $[0, T + \epsilon]$. Hence,
$$D_c^\gamma(\eta_\epsilon * [(u - u_0)\theta(t)]) = g_{-\gamma} * (\eta_\epsilon * [(u - u_0)\theta(t)]).$$

It follows that

(2.24) $$\eta_\epsilon * [(u - u_0)\theta(t)] = g_\gamma * D_c^\gamma(\eta_\epsilon * [(u - u_0)\theta(t)]).$$

For any $\varphi \in C_c^\infty(-\infty, T)$, there is $\epsilon_0 > 0$ such that for all $\epsilon \in (0, \epsilon_0)$, we have $\tilde{\eta}_\epsilon * \varphi \in C_c^\infty(-\infty, T)$ and the following holds

$$\langle D_c^\gamma(\eta_\epsilon * (u - u_0)\theta), \varphi \rangle = \langle (u - u_0)\theta, \tilde{\eta}_\epsilon * \tilde{g}_{-\gamma} * \varphi \rangle$$
$$= \langle (u - u_0)\theta, \tilde{D}_{c;T}^\gamma(\tilde{\eta}_\epsilon * \varphi) \rangle = \langle D_c^\gamma u, \tilde{\eta}_\epsilon * \varphi \rangle = \langle \eta_\epsilon * D_c^\gamma u, \varphi \rangle.$$

It follows that in $\mathscr{D}'$ it holds that

$$\lim_{\epsilon \to 0} D_c^\gamma(\eta_\epsilon * (u - u_0)\theta) = \lim_{\epsilon \to 0} \eta_\epsilon * D_c^\gamma u.$$

Using (2.24), we find in $\mathscr{D}'$ that

$$\lim_{\epsilon \to 0} \eta_\epsilon * (u - u_0)\theta = \lim_{\epsilon \to 0} g_\gamma * D_c^\gamma(\eta_\epsilon * [(u - u_0)\theta(t)]) = \lim_{\epsilon \to 0} g_\gamma * \eta_\epsilon * D_c^\gamma u.$$

The claim then follows. □

In general $D_c^\gamma u \in \mathscr{D}'$ is a functional from $C_c^\infty((-\infty, T); \mathbb{R})$ to $B$. We say $D_c^\gamma u \in L^1_{\text{loc}}([0, T); B)$ if there exists a function $f \in L^1_{\text{loc}}([0, T); B)$ such that for any $\varphi \in C_c^\infty((-\infty, T); \mathbb{R})$, we have

$$(2.25) \qquad \langle D_c^\gamma u, \varphi \rangle = \int_0^T f(t)\varphi(t)\, dt.$$

In this case, we will identify $D_c^\gamma u$ with $f$, while identifying the pairing between $C_c^\infty((-\infty, T); \mathbb{R})$ and $\mathscr{D}'$ with the integral (2.25). With this notion, we have the following observation, for which we omit the proof:

PROPOSITION 2.15. *Let $\gamma \in (0, 1)$. If $f := D_c^\gamma u \in L^1_{\text{loc}}([0, T); B)$, then*

$$u(t) = J_\gamma(f)(t) = u_0 + \frac{1}{\Gamma(\gamma)} \int_0^t (t - s)^{\gamma - 1} f(s)\, ds, \text{ a.e. on } (0, T),$$

*where the integral is understood in the Lebesgue sense.*

COROLLARY 2.16. *If the weak Caputo derivative associated with initial value $u_0$ satisfies $D_c^\gamma u \in L^{1/\gamma}_{\text{loc}}([0, T), B)$, then $u(0+) = u_0$ in the sense of Definition 2.1.*

*Proof.* We can estimate directly that as $t \to 0+$,

$$\frac{1}{t} \int_0^t \|u(t) - u_0\|_B dt \leq \frac{1}{t\Gamma(\gamma)} \int_0^t \int_0^\tau (\tau - s)^{\gamma - 1} \|D_c^\gamma u\|_B ds d\tau$$
$$\leq \frac{1}{t\Gamma(\gamma)} \int_0^t \|D_c^\gamma\|_B(s) \int_s^t (\tau - s)^{\gamma - 1} d\tau ds$$
$$= \frac{1}{t\Gamma(1 + \gamma)} \int_0^t (t - s)^\gamma \|D_c^\gamma u\|_B(s)\, ds$$
$$\leq \frac{1}{t\Gamma(1 + \gamma)} \|(t - s)^\gamma\|_{L^{1/(1-\gamma)}(0,t)} \|D_c^\gamma u\|_{L^{1/\gamma}((0,t);B)}$$
$$\leq \frac{1}{\Gamma(1 + \gamma)} \|D_c^\gamma u\|_{L^{1/\gamma}((0,t);B)} \to 0.$$

The last term goes to zero because $\|D_c^\gamma u\|_B^{1/\gamma}$ is integrable on $[0, T - \delta]$ for some $\delta > 0$. □

The following property verifies that our definition agrees with the traditional Caputo derivative if the function has enough regularity.

LEMMA 2.17. *Let $\gamma \in (0,1)$. If $u$ is absolutely continuous on $(0,T)$, then $D_c^\gamma u \in L^1_{\text{loc}}([0,T); B)$ and*

$$(2.26) \qquad D_c^\gamma u(t) = \frac{1}{\Gamma(1-\gamma)} \int_0^t \frac{\dot{u}(s)}{(t-s)^\gamma} ds, \ t \in [0,T).$$

*Proof.* We just need to check that the expression given here satisfies the definition. Since $u$ is absolutely continuous, then $\dot{u} \in L^1((0,T); B)$. Then, by Young's inequality,

$$f(t) := \frac{1}{\Gamma(1-\gamma)} \int_0^t \frac{\dot{u}(s)}{(t-s)^\gamma} ds \in L^1((0,T); B).$$

We compute directly that

$$\frac{1}{\Gamma(1-\gamma)} \int_0^T \varphi(t) \int_0^t \frac{\dot{u}(s)}{(t-s)^\gamma} ds\, dt = \frac{1}{\Gamma(1-\gamma)} \int_0^T \dot{u}(s) \int_s^T \frac{\varphi(t)}{(t-s)^\gamma} dt\, ds$$

$$= \frac{-1}{\Gamma(1-\gamma)} \int_0^T (u(t) - u(0+)) \frac{d}{ds} \int_s^T \frac{\varphi(t)}{(t-s)^\gamma} dt\, ds.$$

Recall that $\varphi \in C_c^\infty((-\infty, T); \mathbb{R})$, we can do integration by parts. Using again that $\varphi(t)$ vanishes at $T$,

$$\frac{d}{ds} \int_s^T \frac{\varphi(t)}{(t-s)^\gamma} dt = \int_s^T \frac{\dot{\varphi}(t)}{(t-s)^\gamma} dt.$$

This verifies that $f$ is the Caputo derivative. $\square$

The following is similar as Lemma 2.4. We omit the proof here.

PROPOSITION 2.18. *Let $\gamma \in (0,1)$. If the mapping $u : [0,T) \to B$ satisfies*

$$u \in C^1((0,T); B) \cap C([0,T); B),$$

*and $u \mapsto E(u) \in \mathbb{R}$ is a $C^1$ convex functional on $B$, then*

$$(2.27) \qquad D_c^\gamma u(t) = \frac{1}{\Gamma(1-\gamma)} \left( \frac{u(t) - u(0)}{t^\gamma} + \gamma \int_0^t \frac{u(t) - u(s)}{(t-s)^{\gamma+1}} ds \right)$$

*and*

$$(2.28) \qquad D_c^\gamma E(u(t)) \leq \langle D_u E(u), D_c^\gamma u \rangle,$$

*where $D_u E(\cdot) : B \to B'$ is the Fréchet differential and $\langle \cdot, \cdot \rangle$ is understood as the dual pairing between $B'$ and $B$.*

**3. Functions with weak Caputo derivatives in $L^p$ and Hölder spaces.** By Proposition 2.15, the functions can be recovered by the fractional integrals of its weak Caputo derivatives. Hence, we can study the time regularity of the functions by studying the regularity improvement of fractional integrals.

The following were proved by Hardy and Littlewood for fractional integrals [39]:

LEMMA 3.1. *Let $B$ be a Banach space and $T > 0$. Suppose*

$$(3.1) \qquad f := D_c^\gamma u \in L^1_{\text{loc}}([0,T); B).$$

(i) If $f \in L^1((0,T); B)$, then $\forall \epsilon \in (0, \frac{\gamma}{1-\gamma}]$,

(3.2) $$\|u - u_0\|_{L^{\frac{1}{1-\gamma} - \epsilon}((0,T);B)} \leq K\|f\|_{L^1((0,T);B)}.$$

(ii) If $f \in L^p((0,T); B)$ for some $p \in (1, 1/\gamma)$, then

(3.3) $$\|u - u_0\|_{L^{\frac{p}{1-p\gamma}}(0,T;B)} \leq K\|f\|_{L^p((0,T);B)}.$$

(iii) If $f \in L^p((0,T); B)$ for some $p > 1/\gamma$, then $u$ continuous on $[0,T]$ such that

(3.4) $$\|u(t+h) - u(t)\|_B \leq Ch^{\gamma - 1/p},$$

for $0 \leq t < t + h \leq T$ and $C$ is independent of $t$.

For the convenience of the readers, the first result of Lemma 3.1 is Section 3.5 (iii) in [39], the second result is Theorem 4 and the third result is Theorem 12 in [39].

We now focus on the regularity of functions with weak Caputo derivatives in Hölder spaces $C^{m,\beta}(U)$, $\beta > 0$ (see [35, sect. 5.1], [7, Chap. 1]). Recall that $f \in C^{m,\beta}(U), \beta \in (0,1]$ means that $f \in C^m(U)$ and $v = f^{(m)}$ satisfies

$$\sup_{x,y \in U, x \neq y} \frac{|v(x) - v(y)|}{|x-y|^\beta} < \infty.$$

If $\beta = 0$, we define $C^{m,\beta} := C^m$, where $C^m$ means the set of $m$-th order continuously differentiable functions. $C^{m,1}$ means $f^{(m)}$ is Lipschitz continuous and clearly $C^{m+1} \subset C^{m,1}$.

It turns out that $C^{m,\beta}$ is sometimes not convenient to use if $\beta = 1$. We introduce the Hölder space $C^{m,\beta;k}, k > 0$ [7, Def. 1.7], which means $f \in C^m$, and $v = f^{(m)}$ satisfies

$$|v(x) - v(y)| \leq C|x-y|^\beta |\ln|x-y||^k, \quad |x-y| < 1/2.$$

Note that we use different notations from [7] to distinguish with the Sobolev spaces $H^\alpha$.

From Lemma 3.1, we can easily infer that if $f \in C([0,T]; B)$, then

$$u(t) = u_0 + \frac{1}{\Gamma(\gamma)} \int_0^t (t-s)^{\gamma-1} f(s)\,ds$$

is $(\gamma - \epsilon)$-Hölder continuous for any $\epsilon > 0$. However, the order of Hölder continuity generally cannot be improved higher than $\gamma$. For example, if $f = 1$ which is smooth, then $u(t) = u_0 + C_1 t^\gamma$ which is only $\gamma$-Hölder continuous at $t = 0$, though smooth for $t > 0$. Actually, this observation is quite general. We have

LEMMA 3.2 ([7], Theorem 3.1). *Suppose $f \in C^{0,\beta}([0,T]; B)$ for a Banach space $B$ and $0 \leq \beta \leq 1$. Let $\gamma \in (0,1)$ and*

$$u(t) = u_0 + \frac{1}{\Gamma(\gamma)} \int_0^t (t-s)^{\gamma-1} f(s)\,ds.$$

*Then,*

(3.5) $$u(t) = u_0 + \frac{f(0)}{\Gamma(1+\gamma)} t^\gamma + \psi(t),$$

*where*

$$\psi \in \begin{cases} C^{0,\beta+\gamma}([0,T];B), & \beta+\gamma < 1, \\ C^{1,\beta+\gamma-1}([0,T];B), & \beta+\gamma > 1, \\ C^{0,1;1}([0,T];B), & \beta+\gamma = 1. \end{cases} \tag{3.6}$$

We have the following results about the regularity improvement:

PROPOSITION 3.3. *Let $B$ be a Banach space and $T > 0$ and $\gamma \in (0,1)$. Suppose $u \in L^1_{\text{loc}}((0,T);B)$ and $f = D^\gamma_c u$ with an assigned initial value $u_0 \in B$.*
  (i) *If $f \in L^\infty((0,T);B)$, then $u$ is $(\gamma - \epsilon)$-Hölder continuous for any $\epsilon \in (0,\gamma)$. If $f$ is continuous, then $u$ is $\gamma$-Hölder continuous.*
  (ii) *If there exists $\delta > 0$, such that $f \in C^{m,\beta}([\delta/4, T]; B)$, with $\beta \in [0,1]$, then*

$$u \in \begin{cases} C^{m,\beta+\gamma}([\delta,T];B), & \beta+\gamma < 1, \\ C^{m+1,\beta+\gamma-1}([\delta,T];B), & \beta+\gamma > 1, \\ C^{m,1;1}([\delta,T];B), & \beta+\gamma = 1. \end{cases} \tag{3.7}$$

*The claims are not true in general if $\delta = 0$.*
  (iii) *If there exists $\delta > 0$, such that $f \in H^s((\delta/4, T); B)$ ($H^s$ denotes the Sobolev space $W^{s,2}$ in time), then*

$$u \in H^{s+\gamma}((\delta, T); B).$$

*The claim is not true in general if $\delta = 0$.*

*Proof.* (i) is the result of Lemma 3.1 (iii) and Lemma 3.2.

For (ii) and (iii), we do the decomposition

$$f = f_1 + f_2$$

so that $\operatorname{supp} f_1 \subset [0, 3\delta/4]$ while $\operatorname{supp} f_2 \subset [\delta/2, T]$ and $f_2$ is again in $C^{m,\beta}([\delta/4, T]; B)$ or in $H^s([\delta/4, T]; B)$. This is doable, for example, by setting $f_i = f\zeta_i$, $i = 1,2$ where $\zeta_i$ are smooth functions such that $\operatorname{supp} \zeta_1 \subset (-\infty, 3\delta/4]$, $\operatorname{supp} \zeta_2 \subset [\delta/2, +\infty)$, and $\zeta_1 + \zeta_2 = 1$. Then, we have

$$u(t) = \left( u_0 + \frac{1}{\Gamma(\gamma)} \int_0^{3\delta/4} (t-s)^{\gamma-1} f_1(s) ds \right) + \frac{1}{\Gamma(\gamma)} \int_0^t (t-s)^{\gamma-1} f_2(s) ds$$
$$=: u_1(t) + u_2(t).$$

The first term $u_1$ is a smooth function on $[\delta, T]$. $u_2$ is treated as follows:

For (ii). we can easily check that $u_2 \in C^m$ and $v = u_2^{(m)}$ satisfies

$$v(t) = \frac{1}{\Gamma(\gamma)} \int_0^t (t-s)^{\gamma-1} f_2^{(m)}(s) \, ds.$$

The claim then follows from Lemma 3.2.

For (iii), the claim follows from [19, Theorem 2.18]. □

Now, we move onto the time shift estimates that are useful for our compactness theorems. We first define the shift operator $\tau_h$ as

$$\tau_h u(t) := u(t+h). \tag{3.8}$$

We have the following claim:

PROPOSITION 3.4. *Fix $T > 0$. Let $B$ be a Banach space and $\gamma \in (0,1)$. Suppose $u \in L^1_{\mathrm{loc}}((0,T); B)$ has a weak Caputo derivative $D^\gamma_c u \in L^p((0,T); B)$ associated with initial value $u_0 \in B$. If $p\gamma \geq 1$, we set $r_0 = \infty$ and if $p\gamma < 1$, we set $r_0 = p/(1-p\gamma)$. Then, there exists $C > 0$ independent of $h$ and $u$ such that*

$$(3.9) \qquad \|\tau_h u - u\|_{L^r((0,T-h);B)} \leq \begin{cases} Ch^{\gamma + \frac{1}{r} - \frac{1}{p}} \|D^\gamma_c u\|_{L^p((0,T);B)}, & r \in [p, r_0), \\ Ch^\gamma \|D^\gamma_c u\|_{L^p((0,T);B)}, & r \in [1, p]. \end{cases}$$

*Proof.* To be convenient, we denote

$$f := D^\gamma_c u \in L^p((0,T); B).$$

By Proposition 2.15, $u(t) = u(0) + \frac{1}{\Gamma(\gamma)} \int_0^t (t-s)^{\gamma-1} f(s)\, ds$.

Denote

$$(3.10) \qquad \begin{aligned} K_1(s,t;h) &:= (t+h-s)^{\gamma-1}, \\ K_2(s,t;h) &:= (t-s)^{\gamma-1} - (t+h-s)^{\gamma-1}. \end{aligned}$$

We then have

$$(3.11) \qquad \tau_h u(t) - u(t) = \frac{1}{\Gamma(\gamma)} \left( \int_t^{t+h} K_1(s,t;h) f(s)\, ds + \int_0^t K_2(s,t;h) f(s)\, ds \right),$$

so that

$$(3.12) \qquad \int_0^{T-h} \|\tau_h u - u\|_B^r\, dt \leq \frac{2^r}{(\Gamma(\gamma))^r} \Big( \int_0^{T-h} \Big( \int_t^{t+h} K_1 \|f\|_B(s)\, ds \Big)^r dt \\ + \int_0^{T-h} \Big( \int_0^t K_2 \|f\|_B(s)\, ds \Big)^r dt \Big).$$

**Case 1:** $r \geq p$ and $\frac{1}{r} > \frac{1}{p} - \gamma$

We denote $I_1 = (t, t+h)$ and $I_2 = (0, t)$. Let $1/r + 1 = 1/q + 1/p$, and we apply Hölder inequality for $i = 1, 2$:

$$(3.13) \qquad \int_{I_i} K_i \|f\|_B(s)\, ds \leq \Big( \int_{I_i} K_i^q \|f\|_B^p\, ds \Big)^{\frac{1}{r}} \Big( \int_{I_i} K_i^q\, ds \Big)^{\frac{r-q}{qr}} \Big( \int_{I_i} \|f\|_B^p\, ds \Big)^{\frac{r-p}{pr}}.$$

We have

$$(3.14) \qquad \Big( \int_{I_i} \|f\|_B^p\, ds \Big)^{\frac{r-p}{pr}} \leq \|f\|_{L^p((0,T);B)}^{1-p/r}.$$

Direct computation shows

$$(3.15) \qquad \int_t^{t+h} K_1^q\, ds = \frac{1}{q(\gamma-1)+1} h^{q(\gamma-1)+1}.$$

Note that for $q \geq 1$, $a \geq 0, b \geq 0$, we have $(a+b)^q \geq a^q + b^q$. Hence,

$$K_2^q \leq (t-s)^{q(\gamma-1)} - (t+h-s)^{q(\gamma-1)}.$$

Since $q(\gamma - 1) + 1 > 0$, we find

(3.16)
$$\int_0^t K_2^q\, ds = \frac{1}{q(\gamma-1)+1}(t^{q(\gamma-1)+1} - (t+h)^{q(\gamma-1)+1} + h^{q(\gamma-1)+1}) \leq Ch^{q(\gamma-1)+1}.$$

Therefore, combining (3.12)-(3.16), we have

$$\int_0^{T-h} \|\tau_h u - u\|_B^r\, dt \leq Ch^{(q(\gamma-1)+1)\frac{r-q}{q}}\Big(\int_0^T ds\|f\|_B^p(s)\int_{0\wedge s-h}^s K_1^q\, dt$$
$$+ \int_0^{T-h} ds\|f\|_B^p(s)\int_s^{T-h} K_2^q\, dt\Big).$$

Direct computation shows $\int_{0\wedge s-h}^s K_1^q\, dt \leq \frac{1}{q(\gamma-1)+1}h^{q(\gamma-1)+1}$ while

$$\int_s^{T-h} K_2^q\, dt \leq \int_s^{T-h}(t-s)^{q(\gamma-1)}\, dt - \int_s^{T-h}(t-s+h)^{q(\gamma-1)}\, dt$$
$$\leq \frac{1}{q(\gamma-1)+1}h^{q(\gamma-1)+1}.$$

Hence,

$$\int_0^{T-h}\|\tau_h u - u\|_B^r\, dt \leq Ch^{(q(\gamma-1)+1)\frac{r}{q}}\|f\|_{L^p((0,T);B)}^r.$$

In other words,

(3.17) $$\|\tau_h u - u\|_{L^r((0,T);B)} \leq Ch^{\gamma+\frac{1}{r}-\frac{1}{p}}\|D_c^\gamma u\|_{L^p((0,T);B)}.$$

**Case 2:** $r < p$

We first note that we have for $r = p$:

$$\|\tau_h u - u\|_{L^p((0,T);B)} \leq Ch^\gamma\|D_c^\gamma u\|_{L^p((0,T);B)}$$

by Case 1.

Then, by Hölder's inequality

(3.18)
$$\|\tau_h u - u\|_{L^r((0,T);B)} \leq \|1\|_{L^{r*(p/(p-r))}((0,T);B)}\|\tau_h u - u\|_{L^{r*(p/r)}((0,T);B)}$$
$$= T^{1/r-1/p}\|\tau_h u - u\|_{L^p((0,T);B)}.$$

This finishes the proof. □

PROPOSITION 3.5. *Suppose $Y$ is a reflexive Banach space, $\gamma \in (0,1)$ and $T > 0$. Assume $u_n \to u$ in $L^p((0,T);Y), p \geq 1$. If there is an assignment of initial values $u_{0,n}$ for $u_n$ such that the weak Caputo derivatives $D_c^\gamma u_n$ are bounded in $L^r((0,T);Y)$ ($r \in [1,\infty)$), then*

(i) *There this a subsequence such that $u_{0,n}$ converges weakly to some value $u_0 \in Y$.*
(ii) *If $r > 1$, there exists a subsequence such that $D_c^\gamma u_{n_k}$ converges weakly to $f$ and $u_{0,n_k}$ converges weakly to $u_0$. Moreover, $f$ is the Caputo derivative of $u$ with initial value $u_0$ so that*

$$u(t) = u_0 + \frac{1}{\Gamma(\gamma)}\int_0^t (t-s)^{\gamma-1}f(s)\, ds.$$

*Further, if $r \geq 1/\gamma$, then, $u(0+) = u_0$ in $Y$ under the sense of Definition 2.1.*

*Proof.* Let $f_n = D_c^\gamma u_n$.

(i). By Lemma 3.1, $u_n(t) - u_{0,n}$ is bounded in $L^{r_1}((0,T);Y)$ where $r_1 \in [1, \frac{r}{1-r\gamma} - \epsilon)$ if $r < 1/\gamma$ or $r_1 \in [1, \infty)$ if $r > 1/\gamma$. Take $p_1 = \min(r_1, p)$. Then, $u_n(t) - u_{0,n}$ is bounded in $L^{p_1}(0,T;Y)$. Since $u_n$ converges in $L^p$ and thus in $L^{p_1}$, then $u_{0,n}$ is bounded in $L^{p_1}((0,T);Y)$. Hence, $u_{0,n}$ is actually bounded in $Y$. Since $Y$ is reflexive, there is a subsequence $u_{0,n_k}$ that converges weakly to $u_0$ in $Y$.

(ii). We can take a subsequence such that both $u_{0,n_k}$ converges weakly to $u_0$ and $f_{n_k} := D_c^\gamma u_{n_k}$ converges to $f$ weakly in $L^r((0,T);Y)$ since $r > 1$. Take $\varphi \in C_c^\infty[0,T)$ and $w \in Y'$. We have

$$(3.19) \qquad \int_0^T \tilde{D}_{c;T}^\gamma(t) \varphi (u_{n_k}(t) - u_{0,n_k}) \, dt = \int_0^T \varphi(t) f_{n_k}(t) \, dt.$$

Since $w\varphi, w\tilde{D}_{c;T}^\gamma \varphi \in L^{r/(r-1)}(0,T;Y')$, we have by (3.19) that

$$\langle w\tilde{D}_{c;T}^\gamma \varphi, u_{n_k}(t) - u_{0,n_k} \rangle = \langle w\varphi, f_{n_k} \rangle,$$

where the pairing is between $L^{r/(r-1)}((0,T);Y')$ and $L^r((0,T);Y)$. Taking the limit $k \to \infty$ and using the weak convergence, we have

$$\left\langle w\tilde{D}_{c;T}^\gamma \varphi, u - u_0 \right\rangle = \langle w\varphi, f \rangle.$$

Since $w$ is arbitrary and $f \in L^r(0,T;Y)$, it must hold that

$$(3.20) \qquad \int_0^T \tilde{D}_{c;T}^\gamma \varphi (u(t) - u_0) \, dt = \int_0^T \varphi f \, dt.$$

Hence, $f$ is the weak Caputo derivative of $u$ with initial value $u_0$. By Proposition 2.15, we have

$$u(t) = u_0 + (g_\gamma * f)(t) = u_0 + \frac{1}{\Gamma(\gamma)} \int_0^t (t-s)^{\gamma-1} f(s) \, ds.$$

The last claim follows from Corollary 2.16. □

**4. Compactness criteria for time fractional PDEs.** For linear evolutionary equations, establishing existence of weak solutions is relatively easy. Indeed, one only needs the weak compactness, which is guaranteed by boundedness in reflexive spaces. For integro-differential equations, one may refer to [32]. However, for nonlinear evolutionary equations, strong compactness criteria like the Aubin–Lions lemma are often needed. In this section, we present and prove some strong compactness criteria which may not be sharp, but are useful for time fractional PDEs.

THEOREM 4.1. *Let $T > 0, \gamma \in (0,1)$ and $p \in [1, \infty)$. Let $M, B, Y$ be Banach spaces. $M \hookrightarrow B$ compactly and $B \hookrightarrow Y$ continuously. Suppose $W \subset L^1_{\text{loc}}((0,T);M)$ satisfies:*

*(i) There exists $C_1 > 0$ such that $\forall u \in W$,*

$$(4.1) \qquad \sup_{t \in (0,T)} J_\gamma(\|u\|_M^p) = \sup_{t \in (0,T)} \frac{1}{\Gamma(\gamma)} \int_0^t (t-s)^{\gamma-1} \|u\|_M^p(s) \, ds \leq C_1.$$

*(ii) There exist $r \in (\frac{p}{1+p\gamma}, \infty) \cap [1, \infty)$ and $C_3 > 0$ such that for any $u \in W$, there is an assignment of initial value $u_0$ to make the weak Caputo derivative satisfy:*

$$(4.2) \qquad \|D_c^\gamma u\|_{L^r((0,T);Y)} \leq C_3.$$

Then, $W$ is relatively compact in $L^p((0, T); B)$.

Another compactness theorem is as following:

THEOREM 4.2. *Let $T > 0, \gamma \in (0, 1)$ and $p \in [1, \infty)$. Let $M, B, Y$ be Banach spaces. $M \hookrightarrow B$ compactly and $B \hookrightarrow Y$ continuously. Suppose $W \subset L^1_{\text{loc}}((0, T); M)$ satisfies:*
*(i). There exists $r_1 \in [1, \infty)$ and $C_1 > 0$ such that $\forall u \in W$,*

$$(4.3) \qquad \sup_{t \in (0,T)} J_\gamma(\|u\|_M^{r_1}) = \sup_{t \in (0,T)} \frac{1}{\Gamma(\gamma)} \int_0^t (t-s)^{\gamma-1} \|u\|_M^{r_1}(s)\, ds \leq C_1.$$

*(ii). There exists $p_1 \in (p, \infty]$ such that $W$ is bounded in $L^{p_1}((0, T); B)$.*
*(iii). There exist $r_2 \in [1, \infty)$ and $C_2 > 0$ such that $\forall u \in W$, there is an assignment of initial value $u_0$ so that the weak Caputo derivative satisfies:*

$$(4.4) \qquad \|D_c^\gamma u\|_{L^{r_2}((0,T);Y)} \leq C_2,$$

*Then, $W$ is relatively compact in $L^p((0, T); B)$.*

To prove the theorems, we need several preliminary results.

**4.1. Bounded fractional integrals.** Regarding the fractional integral, we find it convenient to define $\|\cdot\|_{L^p_\gamma}$ for $\gamma \in (0, 1)$ and $p \geq 1$ as

$$(4.5) \qquad \|u\|_{L^p_\gamma(0,T;M)} := \sup_{t \in (0,T)} \left( \left| \int_0^t (t-s)^{\gamma-1} \|u\|_M^p(s)\, ds \right| \right)^{1/p} < \infty.$$

It is easy to verify that the mapping $\|\cdot\|_{L^p_\gamma(0,T;M)}$ satisfies the triangle inequality so that it is a norm. A simple observation is

LEMMA 4.3. *Let $\gamma \in (0, 1)$. If $\|f\|_{L^p_\gamma((0,T);\|\cdot\|_M)} < \infty$, then $f \in L^p((0, T); M)$.*

*Proof.* The result follows from the simple observation:

$$(4.6) \qquad \int_0^T \|f\|_M^p(s)\, ds \leq T^{1-\gamma} \int_0^T (T-s)^{\gamma-1} \|f\|_M^p\, ds. \qquad \square$$

It is tempting to prove $f \in L^{pr}((0, T); M)$ for some $r > 1$ in Lemma 4.3 by performing more careful estimates. Indeed, if $\{t : \|f(t)\|_M \geq z\}$ is a single interval, we can improve the results, but we also have evidence that the improvement may not be possible in some cases. See Claim 1 for the cases when we can improve and Claim 2 for the evidence below.

**Claim 1.** *If $A_z = \{t \in [0, T] : \|f(t)\|_M \geq z\}$ is a single interval for any $z > 0$ (for example $\|f(\cdot)\|_M$ is monotone) and $\|f\|_{L^p_\gamma((0,T);\|\cdot\|_M)} < \infty$, then $f \in L^{pr}((0, T); M)$ for any $r \in [1, 1/\gamma)$.*

*Proof.* Let $g = \|f(\cdot)\|_M^p$. The fact $\|f\|_{L^p_\gamma((0,T);\|\cdot\|_M)} < \infty$ implies that

$$(4.7) \qquad \sup_{0 \leq t \leq T} \int_0^t (t-s)^{\gamma-1} g(s)\, ds < \infty.$$

The problem is then reduced to showing that $g \in L^r(0, T)$ with $r \in [1, 1/\gamma)$.

For the $L^r(0, T)$ norm of $g$, we have

$$(4.8) \qquad \|g\|_{L^r(0,T)}^r \sim \int_0^\infty z^{r-1} \lambda(z)\, dz,$$

where $\lambda(z) = |\{t : |g(t)| \geq z\}|$. Note that $\{t : |g(t)| \geq z\} = A_{z^{1/p}}$ is an interval. Let the two endpoints be $a_z, b_z$. Then, (4.7) implies that

$$z \int_{a_z}^{b_z} (b_z - s)^{\gamma-1} ds \leq \sup_{0 \leq t \leq T} \int_0^t (t-s)^{\gamma-1} g(s) \, ds \leq C.$$

Then we have $z\lambda(z)^\gamma \leq C_1$ and thus $\lambda(z) \lesssim Cz^{-1/\gamma}$. Hence, (4.8) converges for $r \in [1, 1/\gamma)$. The claim then follows. □

Combining with the results in Lemma 3.1 (ii), the upper bound $1/\gamma$ of $r$ seems to be optimal. One question is whether we can improve $f \in L^{pr}$ for general data. Unfortunately, we believe this is not true.

Recall that a Borel measure is said to be Ahlfors-regular of degree (dimension) $\alpha \in (0, 1)$ [40, Def. 2.11] if there exist $C_1 > 0, C_2 > 0$ such that it holds for all $x \in \operatorname{supp} \mu$ that

(4.9) $$C_1 r^\alpha \leq \mu(B(x,r)) \leq C_2 r^\alpha.$$

**Claim 2.** *Suppose $\mu$ is the middle $1/3$ Cantor measure [41] that is Ahlfors-regular of degree (or dimension) $\alpha = \ln 2/\ln 3$. Then, if $\gamma > 1 - \alpha$,*

$$\sup_{t \in [0,1]} \int_0^1 |t-s|^{\gamma-1} d\mu(s) < \infty.$$

*Proof.* We perform the dyadic decomposition of the interval:

$$I_k = [(1-2^{-k})t, (1-2^{-k-1})t) \cup (t + (1-t)2^{-k-1}, t + (1-t)2^{-k}] =: I_{k1} \cup I_{k2}.$$

Clearly, $\cup_{k=0}^\infty I_k = [0,1] \setminus \{t\}$. Since $\mu\{t\} = 0$, it suffices to show that

$$\sum_{k \geq 0} \int_{I_k} |t-s|^{\gamma-1} d\mu(s) < \infty.$$

If $s \in I_{k1}$, we have

$$|t-s|^{\gamma-1} \leq 2^{(k+1)(1-\gamma)} t^{\gamma-1} = |I_{k1}|^{\gamma-1}.$$

If $s \in I_{k2}$, we have

$$|t-s|^{\gamma-1} \leq 2^{(k+1)(1-\gamma)} (1-t)^{\gamma-1} = |I_{k2}|^{\gamma-1}.$$

It follows that

$$\int_{I_k} |t-s|^{\gamma-1} d\mu(s) \leq C_2(|I_{k1}|^{\alpha+\gamma-1} + |I_{k2}|^{\alpha+\gamma-1})$$
$$= C_2(2^{-(k+1)(\alpha+\gamma-1)} t^{\alpha+\gamma-1} + 2^{-(k+1)(\alpha+\gamma-1)} (1-t)^{\alpha+\gamma-1}).$$

It follows that if $\delta = \alpha + \gamma - 1 > 0$,

$$\int_0^1 |t-s|^{\gamma-1} d\mu(s) \leq C_2 \frac{2^{-\delta}}{1-2^{-\delta}} (t^\delta + (1-t)^\delta) \leq 2C_2 \frac{2^{-\delta}}{1-2^{-\delta}}.$$

□

Using the result in Claim 2, we can have some function $g$ that is close to $\mu$ (we mean the measure $g\,dt$ is close to $\mu$) and $\sup_{0\leq t\leq T}\int_0^t(t-s)^{\gamma-1}g(s)\,ds < \infty$, but the $L^r$ ($r>1$) norm can be as large as possible. In fact, since $\mu$ is supported on a Lebesgue measure zero set, it is clear that $\mu*\eta_\epsilon$ (with $\eta_\epsilon$ being a mollifier) is a Lebesgue measurable function but $\sup_{\epsilon>0}\|\mu*\eta_\epsilon\|_{L^r} = \infty$. This essentially forbids any improvement of the result in Lemma 4.3. Furthermore, for an arbitrary degree $\alpha \in (0,1)$, there is a corresponding Cantor measure so $\alpha = \ln 2/\ln 3$ is not really a critical value.

**4.2. Proof of the compactness criteria.** We first recall some classical results for compact sets in $L^p((0,T);B)$. The first is:

LEMMA 4.4 ([36, Theorem 5]). *Suppose $M, B, Y$ are three Banach spaces. $M \hookrightarrow B \hookrightarrow Y$ with the embedding $M \to B$ being compact. Assume that $1 \leq p < \infty$ and $W \subset L^1_{\mathrm{loc}}((0,T);M)$ satisfies:*
*(i). $W$ is bounded in $L^p((0,T);M)$;*
*(ii). $\|\tau_h f - f\|_{L^p((0,T-h);Y)} \to 0$ uniformly for $f \in W$ as $h \to 0$.*
*Then, $W$ is relatively compact in $L^p((0,T);B)$.*

The second one is

LEMMA 4.5 ([36, Lemma 3]). *Let $1 < p_1 \leq \infty$. If $W$ is a bounded set in $L^{p_1}((0,T);B)$ and relatively compact in $L^1_{\mathrm{loc}}((0,T);B)$, then it is relatively compact in $L^p((0,T);B)$ for all $1 \leq p < p_1$.*

With all the preparation made, we are able to prove the compactness criteria now.

*Proof of Theorem 4.1.* Condition (i) implies that $W$ is bounded in $L^p((0,T);B)$ by Lemma 4.3. Consider Condition (ii). Let $r_0$ be the number in Proposition 3.4. If $r < 1/\gamma$,
$$r_0 = \frac{r}{1-r\gamma} > p,$$
since $r > p/(1+p\gamma)$. Otherwise, $r_0 = \infty > p$. Together with the condition $r \geq 1$, Condition (ii) therefore implies that $\|\tau_h u - u\|_{L^p((0,T);B)} \to 0$ uniformly by Proposition 3.4. Hence, the second condition of Lemma 4.4 is verified. By Lemma 4.4, the relative compactness of $W$ in $L^p((0,T);B)$ follows. □

*Proof of Theorem 4.2.* Condition (i) implies that $W$ is bounded in $L^p((0,T);B)$ by Lemma 4.3. By Condition (iii) and Proposition 3.4, $\|\tau_h u - u\|_{L^1((0,T);B)} \to 0$ uniformly. Hence, by Lemma 4.4, $W$ is relatively compact in $L^1((0,T);B)$. Since $W$ is bounded in $L^{p_1}((0,T);B)$ with $p_1 > p$ according to Condition (ii), the relative compactness of $W$ in $L^p((0,T);B)$ follows from Lemma 4.5. □

**5. Time fractional PDE examples.** In this section, we look at two nonlinear fractional PDEs and see how our compactness theorems can be used to give the existence of weak solutions. The first example is a special case of the time fractional compressible Navier–Stokes equations (with constant density) while the second example is the time fractional Keller–Segel equations.

**5.1. Time fractional compressible Navier–Stokes equations.** The famous Navier–Stokes equations (compressible or incompressible) describe the dynamics of Newtonian fluids [42, 43, 44]. For the incompressible case with constant density, the existence and uniqueness of weak solution in 2D have been proved. However, in 3D case, the global weak solutions may not be unique. The existence and uniqueness of global smooth solutions are still open [45]. For compressible cases, one can refer to [46].

In this subsection, we use the compressible Navier–Stokes equations with constant density as a base model and replace the time derivative with the fractional time derivative. We will use our compactness criteria to show the existence of weak solutions for this model problem.

Let

(5.1) $$\Omega \subset \mathbb{R}^d, \ d = 2, \ 3$$

be a bounded open set with smooth boundary. The following special case of time-fractional compressible Navier–Stokes equations (for usual time derivative, the equations are also known as Euler-Poincare equations or multi-dimensional Burgers equations [47]) we consider are given by

$$D_c^\gamma u + u \cdot \nabla u + (\nabla u) \cdot u + (\nabla \cdot u)u = \Delta u, \ x \in \Omega,$$
$$u|_{\partial \Omega} = 0.$$

Here, the tensor $\nabla u$ is given by $(\nabla u)_{ij} = \partial_i u_j$. This can also be formulated as the conservative form

(5.2) $$D_c^\gamma u + \nabla \cdot (u \otimes u) + \frac{1}{2}\nabla(|u|^2) = \Delta u,$$
$$u|_{\partial \Omega} = 0.$$

The tensor product $u \otimes u$ is given by $(u \times u)_{ij} = u_i u_j$.

**5.1.1. Weak formulation.** Motivated by the integration by parts (Lemma 2.8 and Definition 2.9), we define

DEFINITION 5.1. *Let $\gamma \in (0, 1)$. We say*

$$u \in L^\infty((0,T); L^2(\Omega)) \cap L^2((0,T); H_0^1(\Omega))$$

*with*

$$D_c^\gamma u \in L^{q_1}((0,T); H^{-1}(\Omega)), \ q_1 = \min(2, 4/d),$$

*is a weak solution to* (5.2) *with initial data $u_0 \in L^2(\Omega)$, if*

(5.3) $$\int_0^T \int_\Omega (u(x,t) - u_0) \tilde{D}_{c;T}^\gamma \varphi \, dxdt - \int_0^T \int_\Omega \nabla \varphi : u \otimes u \, dxdt$$
$$- \frac{1}{2} \int_0^T \int_\Omega \nabla \cdot \varphi |u|^2 \, dxdt = \int_0^T \int_\Omega u \cdot \Delta \varphi \, dxdt,$$

*for any $\varphi \in C_c^\infty([0,T) \times \Omega; \mathbb{R}^d)$. We say a weak solution is a regular weak solution if $u(0+) = u_0$ under $H^{-1}(\Omega)$ in the sense of Definition 2.1.*

*If $u$ is a function defined on $(0, \infty)$ so that its restriction on any interval $[0, T)$, $T > 0$ is a (regular) weak solution, we say $u$ is a global (regular) weak solution.*

**Remark 5.1.** *Usually, the test functions $\varphi$ are chosen in a suitable Banach space that makes all the integrals meaningful. The smooth test functions, however, are general enough by a density argument.*

**5.1.2. Preliminary a priori estimates.** Note that if we assume that Proposition 2.18 holds for $u$ and note that $\frac{1}{2}\|u\|_2^2$ is a convex functional, we have

$$D_c^\gamma \frac{1}{2}\|u\|_2^2 \leq \langle u, D_c^\gamma u \rangle = -\int_\Omega \nabla \cdot \left(\frac{1}{2}|u|^2 u\right) dx - \int_\Omega |\nabla u|^2 dx = -\int_\Omega |\nabla u|^2 dx.$$

In other words,

$$D_c^\gamma \frac{1}{2}\|u\|_2^2 \leq -\|\nabla u\|_{L^2}^2.$$

We have therefore by Lemma 2.3 that

(5.4) $$\frac{1}{2}\|u(t)\|_2^2 + \frac{1}{\Gamma(\gamma)} \int_0^t (t-s)^{\gamma-1} \|\nabla u\|_{L^2}^2(s) ds \leq \frac{1}{2}\|u_0\|_2^2.$$

Consequently, $u \in L^\infty((0,T); L^2(\Omega)) \cap L^2((0,T); H_0^1(\Omega))$.

Consider that $d = 2, 3$. Let $p_1 = \max\left(2, \frac{4}{4-d}\right)$, and $q_1 = \min(2, 4/d)$ is the conjugate index of $p_1$. Let $\varphi \in L^{p_1}(0,T; H_0^1(\Omega))$.

(5.5) $$|\langle D_c^\gamma u, \varphi \rangle_{x,t}| = \left|\left\langle -\nabla \cdot (uu) - \frac{1}{2}\nabla(|u|^2) + \Delta u, \varphi \right\rangle_{x,t}\right|$$

$$\leq C \int_0^T \|\nabla \varphi |u|^2\|_1 \, dt + \int_0^T \|\nabla \varphi\|_2 \|\nabla u\|_2 \, dt.$$

Using the Gagliardo–Nirenberg inequality $\|u\|_4 \leq C\|u\|_2^{1-d/4} \|Du\|_2^{d/4}$, the first term is estimated as

(5.6) $$\int_0^T \|\nabla \varphi |u|^2\|_1 \, dt \leq \int_0^T \|\nabla \varphi\|_2 \|u\|_4^2 \, dt$$

$$\leq C \left(\int_0^T \|\nabla \varphi\|_2^{4/(4-d)} dt\right)^{(4-d)/4} \left(\int_0^T \|Du\|_2^2 dt\right)^{d/4},$$

where we have used the fact $\|u\|_{L^\infty(0,T;L^2)} \leq \|u_0\|_2$. It is then clear that

(5.7) $$D_c^\gamma u \in L^{q_1}((0,T); H^{-1}(\Omega)).$$

Theorem 4.2 can be used to give the compactness for the approximation sequences if these a priori estimates are preserved for the approximation sequences.

**5.1.3. Existence of weak solutions: a Galerkin method.** With the a priori energy estimates, the existence of weak solutions can be performed by the standard techniques. We first of all state the results as follows.

THEOREM 5.2. *Suppose $u_0 \in L^2(\Omega)$. Then there exists a global weak solution to Equation (5.2) with initial data $u_0$ under Definition 5.1. Further, if $\max(\frac{1}{2}, \frac{d}{4}) \leq \gamma < 1$, there is a global weak solution continuous at $t = 0$ under the $H^{-1}(\Omega)$ norm, and hence a global regular weak solution.*

**Remark 5.2.** To guarantee the continuity of $u$ at $t = 0$ for all $\gamma$, we need $D_c^\gamma u \in L^{q_1}((0,T); H^{-1}(\Omega))$ for all $q_1 \in (1, \infty)$. Claim 2 forbids us to conclude that $\nabla u \in L^r(0,T; L^2)$ for $r > 2$. We cannot improve $q_1$ even if we consider the weaker norm of $D_c^\gamma u$ (for example, $L^{q_1}(0,T; H^\alpha(\Omega))$ norm with $\alpha < -1$) due to $\nabla \varphi |u|^2$ term in (5.6).

We use Galerkin method to prove this. Let $\{w_n\}_{n=1}^{\infty}$ be a basis of both $H_0^1(\Omega)$ and $L^2(\Omega)$, and orthonormal in $L^2(\Omega)$, which as well-known exists (see [35, sect. 6.5]).

We first introduce $P_m$ to be the projection that projects $v$ onto the first $m$ modes. In other words, for any $v \in H_0^1(\Omega)$,

$$v = \sum_{k=1}^{\infty} \alpha_k w_k, \tag{5.8}$$

we define

$$P_m v := \sum_{k=1}^{m} \alpha_k w_k. \tag{5.9}$$

Note that one can first of all assume (5.8) holds in $H_0^1$. Then, it must also hold in $L^2$ since $H_0^1 \subset L^2$. This means that the expansion coefficients of $v$ in $H_0^1$ and $L^2$ are the same. We first note that for a general Banach space $B$ and a Schauder basis $\{w_k\} \subset B$, the so-defined projection $P_m$ is a bounded operator, and one may refer to [48, Page 32]. As a consequence of the uniform boundedness principle (Banach–Steinhaus theorem) [49], one has the following well-known fact:

LEMMA 5.3. *Suppose $\{w_k\}_{k=1}^{\infty}$ is a Schauder basis of a Banach space $B$. Consider the projection operator $\{P_m\}$ as in (5.9). Then, $P_m : B \to B$ is a bounded linear operator and*

$$\sup_{m \geq 1} \|P_m\| < \infty,$$

*where*

$$\|P_m\| := \sup_{v \in B, v \neq 0} \frac{\|P_m v\|_B}{\|v\|_B}. \tag{5.10}$$

Let $u_0 = \sum_{k=1}^{\infty} \alpha^k w_k(x)$ in $L^2(\Omega)$. We pursue the function

$$u_m(t) = \sum_{k=1}^{m} c_m^k(t) w_k \tag{5.11}$$

such that $c_m := (c_m^1, \ldots, c_m^m)$ is continuous in time and $u_m$ satisfies the following equations

$$\langle D_c^\gamma u_m, w_j \rangle + \langle \nabla \cdot (u_m \otimes u_m), w_j \rangle + \frac{1}{2} \langle \nabla |u_m|^2, w_j \rangle = \langle \Delta u_m, w_j \rangle,$$

$$u_m(0) = \sum_{k=1}^{m} c_m^k(0) w_k = \sum_{k=1}^{m} \alpha^k w_k. \tag{5.12}$$

Since $c_m$ is continuous, $D_c^\gamma u_m$ is the Caputo derivative (with natural initial value).

Equation (5.12) can be reduced to the following FODE system for $c_m$

$$D_c^\gamma c_m = F_m(c_m),$$

$$c_m(0) = (\alpha^1, \ldots, \alpha^m), \tag{5.13}$$

where $F_m$ is clearly a quadratic vector-valued function of $c_m$, and hence smooth. By studying the FODE system (5.13), we have LEMMA 5.4. *(i) For any $m \geq 1$, there exists a unique solution $u_m$ of the form (5.11) to (5.12) that is continuous on $[0, \infty)$, satisfying*

(5.14)
$$\|u_m\|_{L^\infty((0,\infty);L^2(\Omega))} \leq \|u_0\|_2, \quad \sup_{0 \leq t < \infty} \int_0^t (t-s)^{\gamma-1} \|\nabla u_m\|_2^2 \, ds \leq \frac{1}{2}\Gamma(\gamma)\|u_0\|_2.$$

*(ii) There exists $u \in L^\infty((0,\infty); L^2(\Omega)) \cap L^2_{\text{loc}}([0,\infty); H_0^1(\Omega))$ and a subsequence $m_k$ such that*
$$u_{m_k} \to u, \text{ in } L^2_{\text{loc}}([0,\infty); L^2(\Omega)).$$

*Further, $u$ has a weak Caputo derivative $D_c^\gamma u \in L^{q_1}_{\text{loc}}([0,\infty); H^{-1}(\Omega))$, where $q_1 = \min\left(2, \frac{4}{d}\right)$.*

*Proof.* (i). By the results for FODE in [19], $c_m(t)$ exists on $[0, T_b^m)$ where either $T_b^m = \infty$ or $T_b^m < \infty$ and $\limsup_{t \to T_b^m -} |c_m| = \infty$ where $|c_m| = \sqrt{\sum_j (c_m^j)^2}$. Note that the norm for $c_m$ is not important because any norms are equivalent for finite dimensional vectors. Further, since $F_m$ is quadratic, by [20, Lemma 3.1], $c_m \in C^1(0, \infty) \cap C[0, \infty)$ and consequently,
$$u_m \in C^1((0, T_b^m); H_0^1(\Omega)) \cap C([0, T_b^m); H_0^1(\Omega)).$$

By Proposition 2.18, we have
$$D_c^\gamma \left(\frac{1}{2}\|u_m\|_2^2\right)(t) \leq \langle u_m, D_c^\gamma u_m \rangle.$$

Since $u_m = \sum_{k=1}^m c_m^k(t) w_k$, using (5.12):

(5.15)
$$\langle u_m, D_c^\gamma u_m \rangle + \int_\Omega u_m \cdot \nabla \cdot (u_m \otimes u_m) \, dx + \frac{1}{2} \int_\Omega u_m \cdot \nabla |u_m|^2 dx = -\int_\Omega |\nabla u_m|^2 dx.$$

Hence, we have
$$D_c^\gamma \left(\frac{1}{2}\|u_m\|_2^2\right)(t) \leq \langle u_m, D_c^\gamma u_m \rangle = -\|\nabla u_m\|_2^2.$$

This implies that
$$\|u_m\|_2^2 + \frac{2}{\Gamma(\gamma)} \int_0^t (t-s)^{\gamma-1} \|\nabla u_m(s)\|_2^2 \, ds \leq \|u_0\|_2^2.$$

Consequently, we find that $T_b^m = \infty$. The first claim follows.

(ii).
Take a test function $v \in L^{p_1}((0,T); H_0^1)$ ($p_1 = \max(2, \frac{4}{4-d})$) with
$$\|v\|_{L^{p_1}((0,T);H_0^1)} \leq 1.$$

Denote

(5.16)
$$v_m := P_m v,$$

where $P_m$ is defined in (5.9). Then, by Lemma 5.3, there exists $C(\Omega, T)$ such that

$$\|v_m\|_{L^{p_1}((0,T);H_0^1)} \leq C.$$

We have

$$\langle D_c^\gamma u_m, v \rangle = \langle D_c^\gamma u_m, v_m \rangle = -\langle \nabla \cdot (u_m \otimes u_m), v_m \rangle - \frac{1}{2}\langle \nabla |u_m|^2, v_m \rangle + \langle \Delta u_m, v_m \rangle.$$

Note that the second equality holds because $v_m \in \text{span}\{w_1, \ldots, w_m\}$.

Using similar tricks as we did in Equations (5.5)-(5.6), we find:

(5.17) $$\langle D_c^\gamma u_m, v \rangle \leq C, \ q_1 = \min\left(2, \frac{4}{d}\right).$$

Hence, $\|D_c^\gamma u_m\|_{L^{q_1}(0,T;H^{-1})} \leq C$ for all $m$.

Now, we have

$$\sup_{0 \leq t \leq T} J_\gamma(\|\nabla u_m\|_2^2) \leq C, \ u_m \in L^\infty(0,T;L^2(\Omega)), \ \|D_c^\gamma u_m\|_{L^{q_1}(0,T;H^{-1})} \leq C.$$

By Theorem 4.2, there is a subsequence $\{u_{m_k}\}$ that converges in $L^p((0,T);L^2(\Omega))$ for any $p \in [1,\infty)$. In particular, we choose $p = 2$.

According to Proposition 3.5, $u$ has a weak Caputo derivative with initial value $u_0$ such that

$$D_c^\gamma u \in L^{q_1}((0,T);H^{-1}(\Omega)).$$

By a standard diagonal argument, $u$ is defined on $(0,\infty)$ and $D_c^\gamma u \in L^{q_1}_{\text{loc}}([0,\infty);H^{-1})$, such that

$$u_{m_k} \to u, \text{ in } L^2_{\text{loc}}([0,\infty);L^2(\Omega)).$$

By taking a further subsequence, we can assume that $u_{m_k}$ also converges a.e. to $u$ in $[0,\infty) \times \Omega$. It is easy to see that

$$\int_{t_1}^{t_2} \|u_m\|_2^2 \, dt \leq \|u_0\|_2^2 (t_2 - t_1).$$

According to Fatou's lemma, we find

$$\int_{t_1}^{t_2} \|u\|_2^2 \, dt \leq \|u_0\|_2^2 (t_2 - t_1),$$

for any $t_1 < t_2$. This then implies that $u \in L^\infty((0,\infty);L^2(\Omega))$.

Fix any $T > 0$. Since $u_{m_k}$ is bounded in $L^2((0,T);H_0^1(\Omega))$, it has a further subsequence that converges weakly in $L^2((0,T);H_0^1(\Omega))$. By a standard diagonal argument, there is a subsequence that converges weakly in $L^2_{\text{loc}}([0,\infty);H_0^1(\Omega))$. The limit must be $u$ by pairing with a smooth test function. Hence, $u \in L^2_{\text{loc}}([0,\infty);H_0^1(\Omega))$. □

**Remark 5.3.** *One may want to show that $u \in L^2_\gamma(0,T;H_0^1(\Omega))$ since $u_m \in L^2_\gamma(0,T;H_0^1(\Omega))$ for all $m$. To this end, we may want to prove that the space $L^2_\gamma$ defined in Section 4.1 is reflexive, which is left for future.*

Now, we can prove Theorem 5.2:

*Proof of Theorem 5.2.* By Lemma 5.4, there is a subsequence that converges in $L^p((0,T); L^2(\Omega))$ for any $p \in [1,\infty)$. Let the limit function be $u$.

Now, for any test function $\varphi \in C_c^\infty([0,T) \times \Omega; \mathbb{R}^d)$, we expand

$$\varphi = \sum_{k=1}^\infty \beta_k w_k, \tag{5.18}$$

and we define

$$\varphi_m := \sum_{k=1}^m \beta_k w_k. \tag{5.19}$$

Since $\varphi$ is a smooth function in $t$ that vanishes at $T$, so is $\varphi_m$, and $\tilde{D}_{c;T}^\gamma \varphi_m \to \tilde{D}_{c;T}^\gamma \varphi$ in $L^{p_1}((0,T); H_0^1)$.

We first of all fix $m_0 \geq 1$, and for $m_j \geq m_0$, we have

$$\langle u_{m_j} - u_0, \tilde{D}_{c;T}^\gamma \varphi_{m_0} \rangle_{x,t} = \langle D_c^\gamma u_{m_j}, \varphi_{m_0} \rangle_{x,t}$$
$$= \int_0^T \int_\Omega \nabla \varphi_{m_0} : u_{m_j} \otimes u_{m_j} \, dxdt$$
$$+ \frac{1}{2} \int_0^T \int_\Omega \nabla \cdot \varphi_{m_0} |u_{m_j}|^2 dxdt - \int_0^T \int_\Omega \nabla \varphi_{m_0} : \nabla u_{m_j} dxdt.$$

The first equality here holds by the integration by parts formula while the second one holds because $\varphi_{m_0} \in \mathrm{span}\{w_1, \ldots, w_{m_j}\}$.

According to the convergence proved in Lemma 5.4, taking $j \to \infty$, we have

$$\int_0^T \int_\Omega (u - u_0) \tilde{D}_{c;T}^\gamma \varphi_{m_0} \, dxdt = \int_0^T \int_\Omega \nabla \varphi_{m_0} : u \otimes u \, dxdt \tag{5.20}$$
$$+ \frac{1}{2} \int_0^T \int_\Omega \nabla \cdot \varphi_{m_0} |u|^2 dxdt - \int_0^T \int_\Omega \nabla \varphi_{m_0} : \nabla u \, dxdt.$$

Then, taking $m_0 \to \infty$, by the convergence $\varphi_m \to \varphi$ in $L^p((0,T); H_0^1(\Omega))$ for any $p \in (1,\infty)$, we find that the weak formulation holds.

Further, if $q_1 \geq 1/\gamma$ or $\gamma \geq \max(1/2, d/4)$, by Lemma 5.4 (ii) and Corollary 2.16, it is a regular weak solution. □

**Remark 5.4.** *For the incompressible fractional Navier–Stokes equations*

$$\begin{cases} D_c^\gamma u + u \cdot \nabla u = -\nabla p + \Delta u, \\ \nabla \cdot u = 0, \end{cases} \tag{5.21}$$

*the existence of weak solutions can also be shown. The a priori estimates follow by dotting $u$ and integrating on $x$:*

$$\frac{1}{2} D_c^\gamma \|u\|_{L^2}^2 \leq -\|\nabla u\|_{L^2}^2.$$

*For the Galerkin approximation, we need to find a basis for the divergence free subspace of $H^1$. Consider the projection operator $P_m$ that projects a function into the subspace spanned by the first $m$ basis functions that are divergence free. Then,*

$$D_c^\gamma u_m + P_m(u_m \cdot \nabla u_m) = \Delta u_m.$$

*Using similar techniques, we can show the existence of weak solutions.*

**5.2. Time fractional Keller–Segel equations.** The Keller–Segel equations are a model for chemotaxis of bacteria [50, 51, 52]. This model has attracted a lot of attention due to its mathematical structures. The weak solutions for Keller–Segel equations in 2D have been totally solved in [52]. The discussion of weak solutions of extended models can be found in [53, 54, 55].

As a toy example for our compactness theory, we replace the usual time derivative in the Keller–Segel equations with the Caputo derivatives and consider the following fractional Keller–Segel equations in $\mathbb{R}^2$:

(5.22) $$\begin{cases} D_c^\gamma \rho + \nabla \cdot (\rho \nabla c) = \Delta \rho, \ x \in \mathbb{R}^2, \\ -\Delta c = \rho, \ x \in \mathbb{R}^2. \end{cases}$$

We first of all introduce the definition of weak solutions.

DEFINITION 5.5. *Given $\rho_0 \geq 0$ and $\rho_0 \in L^1(\mathbb{R}^2) \cap L^2(\mathbb{R}^2)$, we say*

$$\rho \in L^\infty(0,T;L^1(\mathbb{R}^2)) \cap L^\infty(0,T;L^2(\mathbb{R}^2)) \cap L^2(0,T;H^1(\mathbb{R}^2))$$

*is a weak solution to the fractional Keller–Segel equation* (5.22) *with initial data $\rho_0$ if*
 (i) $\rho(x,t) \geq 0$.
 (ii) *There exists $q \in (1,2)$ such that $D_c^\gamma \rho \in L^{q_1}((0,T); W^{-2,q}(\mathbb{R}^2))$ for any $q_1 \in (1,\infty)$.*
 (iii) *For any $\varphi \in C_c^\infty([0,T] \times \mathbb{R}^2)$,*

$$\int_0^T \int_{\mathbb{R}^2} (\rho - \rho_0) \tilde{D}_{c;T}^\gamma \varphi \, dxdt - \int_0^T \int_{\mathbb{R}^2} \nabla \varphi \cdot (\nabla(-\Delta)^{-1}\rho) \rho \, dxdt = \int_0^T \int_{\mathbb{R}^2} \rho \Delta \varphi \, dxdt.$$

*We say a weak solution is a regular weak solution if $\rho(0+) = \rho_0$ under $W^{-2,q}$ in the sense of Definition 2.1, where $q$ is given as in (ii).*

*If $\rho$ is a function defined on $(0,\infty)$ so that its restriction on any interval $[0,T)$ ($T > 0$) is a (regular) weak solution, we say $\rho$ is a global (regular) weak solution.*

To study the existence of weak solutions, we first of all investigate the fractional advection diffusion equation

(5.23) $$D_c^\gamma \rho + \nabla \cdot (\rho a(x,t)) = \Delta \rho,$$

with initial condition

$$\rho(x,0) = \rho_0(x).$$

Introduce the Mittag–Leffler function

(5.24) $$E_\gamma(z) := \sum_{n=0}^\infty \frac{z^n}{\Gamma(n\gamma + 1)}$$

and denote $A = -\Delta$ which is a self-joint positive operator. By taking the Laplace transform of the equation, one has the following analogy of Duhamel's principle (though the dynamics is not Markovian) [25, Sections 8-9]:

(5.25) $$\rho(x,t) = E_\gamma(-t^\gamma A)\rho_0 + \gamma \int_0^t \tau^{\gamma-1} E_\gamma'(-\tau^\gamma A)(-\nabla \cdot (\rho a)|_{t-\tau}) \, d\tau.$$

DEFINITION 5.6. *Suppose $X$ is a Banach space in space and time. If $\rho \in X$ satisfies (5.25), then we say $\rho$ is a mild solution of (5.23) in $X$.*

We have the following lemma regarding (5.23) whose proof is in Appendix A:

LEMMA 5.7. *Suppose $a(x,t)$ is smooth and all the derivatives are bounded. Then:*
(i) *If*
$$\rho_0 \in L^1(\mathbb{R}^2) \cap H^\alpha(\mathbb{R}^2),$$
*then $\forall T > 0$, (5.23) has a unique mild solution in $C([0,T]; H^\alpha(\mathbb{R}^2))$.*
(ii) *For the unique mild solution in (i), $\forall T > 0$,*

$$(5.26) \qquad \rho \in C^{0,\gamma}([0,T]; H^\alpha(\mathbb{R}^2)) \cap C^\infty((0,T); H^\alpha(\mathbb{R}^2)).$$

*Moreover, the following holds strongly in $C([0,T]; H^{\alpha-2})$*

$$(5.27) \qquad D_c^\gamma \rho = \frac{1}{\Gamma(1-\gamma)} \int_0^t \frac{\dot\rho(s)}{(t-s)^\gamma} ds = -\nabla \cdot (\rho a(x,t)) + \Delta \rho.$$

(iii) *If $\rho_0 \in H^1(\mathbb{R}^2) \cap L^1(\mathbb{R}^2)$ and $\rho_0 \geq 0$, then $\rho(x,t) \geq 0$, and*

$$(5.28) \qquad \int_{\mathbb{R}^2} \rho \, dx = \int_{\mathbb{R}^2} \rho_0 \, dx.$$

Based on Lemma 5.7, we are motivated to consider the mollified equation. Let $J(x) \in C_c^\infty(\mathbb{R}^2)$ such that $J(x) \geq 0$ and $\int_{\mathbb{R}^2} J(x) dx = 1$. Introduce

$$(5.29) \qquad J_\epsilon = \frac{1}{\epsilon^2} J\left(\frac{x}{\epsilon}\right).$$

The mollified equations read:

$$(5.30) \qquad \begin{cases} D_c^\gamma \rho^\epsilon + \nabla \cdot (\rho^\epsilon \nabla c^\epsilon) = \Delta \rho^\epsilon, \\ -\Delta c^\epsilon = \rho^\epsilon * J_\epsilon, \end{cases}$$

with initial data

$$\rho_0^\epsilon = \rho_0 * J_\epsilon,$$

which has the same $L^1$ norm as $\rho_0$.

Note that the system (5.30) is nonlinear, but the system is similar to the linear problem (5.23) because $\rho^\epsilon * J_\epsilon$ (and hence $\nabla c^\epsilon$) is smooth and bounded, and its derivatives are bounded. Actually, we have

LEMMA 5.8. *Given $\rho_0 \in L^1(\mathbb{R}^2) \cap L^2(\mathbb{R}^2)$ and $\rho_0 \geq 0$, the regularized system (5.30) has a unique global mild solution. Further, this mild solution $\rho^\epsilon$ is a strong solution and $\rho^\epsilon \in C([0,\infty), C^k(\mathbb{R}^2))$ for any $k \geq 0$, and further $\rho^\epsilon \geq 0$.*

*Proof.* The existence and uniqueness of mild solution is similar as we do for the linear problem (see the proof of Lemma 5.7 in Appendix A) and we omit the details. The key idea is to use formula (5.25) and to notice

$$\nabla c^\epsilon = C \frac{x}{|x|^2} * J_\epsilon * \rho^\epsilon,$$

which is smooth and bounded, with derivatives bounded. For similar discussion, one can refer to [56] but our problem is much easier compared with [56] since $\nabla c^\epsilon$ is bounded.

Now one can consider the following problem with initial data $\rho_0$
$$D_c^\gamma v + \nabla \cdot (v\nabla c^\epsilon) = \Delta v.$$

By Lemma 5.7, this problem has a unique mild solution, which is also a strong solution $v \in C([0,\infty), C^k(\mathbb{R}^2))$ for any $k \geq 0$ and $v \geq 0$. However, $\rho^\epsilon$ is also a mild solution, so we must have $v = \rho^\epsilon$, which implies that $\rho^\epsilon$ is also a strong solution to (5.30) with the desired properties. □

We have the following estimates of $\rho^\epsilon$:

LEMMA 5.9. *Suppose $\rho_0 \geq 0$ satisfies that $\rho_0 \in L^1 \cap L^2$ and $M_0 = \|\rho_0\|_1$ is sufficiently small. Then, $\rho^\epsilon \geq 0$ and for any fixed $T > 0$, there are constants $C(q,T) > 0$ and $C(T) > 0$ such that*

(5.31)
$$\|\rho^\epsilon\|_{L^\infty(0,T;L^q)} \leq C(q,T), \forall q \in [1,2],$$
$$\sup_{0\leq t\leq T} \int_0^t (t-s)^{\gamma-1} \|\nabla \rho^\epsilon\|_2^2 \, ds \leq C(T).$$

*Moreover, there exists $q \in (1,2)$ such that $D_c^\gamma \rho^\epsilon$ is uniformly bounded in the space $L^{q_1}([0,T]; W^{-2,q}(\mathbb{R}^2))$ for any $q_1 \in (1,\infty)$.*

*Proof.* By Lemma 5.8, (5.30) has a strong solution $\rho^\epsilon \in C([0,\infty), C^k(\mathbb{R}^2))$ for any $k \geq 0$, and $\rho^\epsilon \geq 0$.

We now perform the estimates of $\rho^\epsilon$. First of all, it is clear that
$$D_c^\gamma \int_{\mathbb{R}^2} \rho^\epsilon dx = 0 \Rightarrow \|\rho^\epsilon\|_1 = \|\rho_0^\epsilon\|_1 = \|\rho_0\|_1.$$

Since $\rho \mapsto \|\rho\|_q^q$ is convex for $q > 1$, by Proposition 2.18,
$$\frac{1}{q} D_c^\gamma \|\rho^\epsilon\|_q^q \leq \langle (\rho^\epsilon)^{q-1}, D_c^\gamma \rho^\epsilon \rangle = \frac{q-1}{q} \|(\rho^\epsilon)^q \rho^\epsilon * J_\epsilon\|_1 - (q-1)\|\nabla(\rho^\epsilon)^{q/2}\|_2^2.$$

Using Hölder's inequality,
$$\|(\rho^\epsilon)^q \rho^\epsilon * J_\epsilon\|_1 \leq \|\rho^\epsilon * J_\epsilon\|_{q+1} \|(\rho^\epsilon)^q\|_{(q+1)/q} \leq \|\rho^\epsilon\|_{q+1}^{q+1}.$$

For $q = 2$, using Gargliardo-Nirenberg inequality,
$$\|\rho^\epsilon\|_3 \leq C \|\nabla \rho^\epsilon\|_2^{2/3} \|\rho^\epsilon\|_1^{1/3}.$$

Hence,
$$\frac{1}{2} D_c^\gamma \|\rho^\epsilon\|_2^2 \leq (C\|\rho^\epsilon\|_1 - 1)\|\nabla \rho^\epsilon\|_2^2 = (C\|\rho_0\|_1 - 1)\|\nabla \rho^\epsilon\|_2^2.$$

If the initial mass $M_0 = \int_{\mathbb{R}^2} \rho_0 \, dx$ is small enough such that
$$CM_0 - 1 < 0,$$
then we have that $\rho^\epsilon$ is uniformly bounded in $L^\infty(0,T;L^2(\mathbb{R}^2)) \cap L^2_{\gamma,\text{loc}}(0,T;H^1(\mathbb{R}^2))$ according to Lemma 2.3.

Since $\rho^\epsilon$ is uniformly bounded in $L^1 \cap L^2$, it is so in $L^p$ for any $p \in [1,2]$. Since $c^\epsilon = (-\Delta)^{-1} \rho^\epsilon$, we have
$$\nabla c^\epsilon = C_1 \frac{x}{|x|^2} * \rho^\epsilon.$$

By the Hardy-Littlewood-Sobolev inequality,

$$\|\nabla c^\epsilon\|_{2p/(2-p)} \leq C_2 \|\rho^\epsilon\|_p, \ 1 < p < 2.$$

Hence, $\nabla c^\epsilon$ is bounded in $L^\infty(0,T; L^r(\mathbb{R}^2))$ for $r \in (2, \infty)$.

We now take test function $\varphi$ with

$$\|\varphi\|_{L^{p_1}(0,T; W_0^{2,p})} \leq 1, \ p > 2, \ p_1 > 1.$$

Then,

$$\langle D_c^\gamma \rho^\epsilon, \varphi \rangle_{x,t} = \langle \rho^\epsilon \nabla c^\epsilon, \nabla \varphi \rangle_{x,t} + \langle \rho^\epsilon, \Delta \varphi \rangle_{x,t}$$

$$\leq \|\rho^\epsilon\|_{L^\infty(0,T;L^2)} \|\nabla c^\epsilon\|_{L^\infty(0,T;L^{2p/(p-2)})} \int_0^T \|\nabla \varphi\|_{L^p} \, dt$$

$$+ \int_0^T \|\Delta \varphi\|_p \|\rho^\epsilon\|_{L^\infty(0,T;L^q)} \, dt.$$

This means

$$\|D_c^\gamma \rho^\epsilon\|_{L^{q_1}(0,T;W^{-2,q})} \leq C(q_1, q, T). \qquad \square$$

The existence of weak solutions is summarized as follows which is a standard consequence of Lemma 5.9 and Theorem 4.1, and we omit the proof.

THEOREM 5.10. *Let $\gamma \in (0,1)$. If $\rho_0 \geq 0$, $\rho_0 \in L^1(\mathbb{R}^2) \cap L^2(\mathbb{R}^2)$ and the initial mass $M_0 = \int_{\mathbb{R}^2} \rho_0 \, dx$ is sufficiently small, then the time-fractional Keller–Segel equation* (5.22) *with initial data $\rho_0$ has a global (non-negative) regular weak solution.*

**Acknowledgements.** The authors would like to thank Xianghong Chen for providing the example in Claim 2. The work of J.-G. Liu is partially supported by KI-Net NSF RNMS11-07444, NSF DMS-1514826 and NSF DMS-1812573.

**Appendix A. Proof of Lemma 5.7.**

*Proof of Lemma 5.7.* (i). Since $E_\gamma(z)$ is an analytic function in the whole $z \in \mathbb{C}$ plane and

$$E_\gamma'(-s) \sim -C_0 s^{-2}, \text{ as } s \to +\infty,$$

we conclude that

$$\sup_{s \in [0,\infty)} E_\gamma'(-s) s^\sigma \leq C, \forall \sigma \leq 2.$$

Consequently,

$$\|E_\gamma'(-\tau^\gamma A) \nabla f\|_{H^\alpha}^2 \leq C \int_{\mathbb{R}^2} E_\gamma'(-\tau^\gamma |k|^2)^2 |k|^2 |\hat{f}_k|^2 (1 + |k|^{2\alpha}) dk$$

$$\leq C \tau^{-\gamma} \int_{\mathbb{R}^2} |\hat{f}_k|^2 (1 + |k|^{2\alpha}) dk = C \tau^{-\gamma} \|f\|_{H^\alpha}^2.$$

We construct the iterative sequence

(A.1)
$$\rho^0(t) = \rho_0, \ \ \rho^n(t) = E_\gamma(-t^\gamma A)\rho_0 + \gamma \int_0^t \tau^{\gamma-1} E_\gamma'(-\tau^\gamma A)(-\nabla \cdot (\rho^{n-1} a)|_{t-\tau}) \, d\tau.$$

We fix $T > 0$. Define $E^n = \rho^n - \rho^{n-1}$. We can compute directly that

$$\|\rho^1\|_{C[0,T;H^\alpha]} \leq \|\rho_0\|_{H^\alpha}\left(1 + C_1\gamma \int_0^t \tau^{\gamma/2-1}d\tau\right) \leq \|\rho_0\|_{H^\alpha}(1 + 2C_1 T^{\gamma/2}).$$

Consequently,
$$\|E^1\|_{C[0,t;H^\alpha]} \leq M, \forall t \in [0,T].$$

The induction formula reads
$$E^n(x,t) = \gamma \int_0^t \tau^{\gamma-1} E'_\gamma(-\tau^\gamma A)(-\nabla \cdot (E^{n-1}a)|_{t-\tau})d\tau.$$

Hence,
$$\|E^n\|_{C([0,t];H^\alpha)} \leq C_1\gamma \sup_{0\leq z\leq t}\int_0^z \tau^{\gamma/2-1}\|E^{n-1}\|_{C([0,z-\tau];H^\alpha)}d\tau$$
$$= C_1\gamma \int_0^t \tau^{\gamma/2-1}\|E^{n-1}\|_{C([0,t-\tau];H^\alpha)}d\tau = C_2 g_{\gamma/2} * \|E^{n-1}\|_{C([0,\cdot];H^\alpha)}.$$

From this induction formula, we have
$$\|E^2\|_{C[0,t;H^\alpha]} \leq C_2 M g_{\gamma/2+1}(t).$$

By induction
$$\|E^n\|_{C[0,t;H^\alpha]} \leq C_2^{n-1} M g_{(n-1)*\gamma/2+1}(t).$$

It follows that
$$\rho = \rho_0 + \sum_{n=1}^\infty E^n$$

converges in $C([0,T]; H^\alpha(\mathbb{R}^2))$ and in other words $\rho^n \to \rho$ in $C([0,T]; H^\alpha(\mathbb{R}^2))$. Hence, $\rho$ is a mild solution.

The uniqueness follows in a similar way. Consider two mild solutions $\rho_1$ and $\rho_2$. Define $\sigma = \rho_1 - \rho_2$. Then,

$$\sigma(x,t) = \gamma \int_0^t \tau^{\gamma-1} E'_\gamma(-\tau^\gamma A)(-\nabla \cdot (\sigma a)|_{t-\tau})d\tau.$$

Then,
$$\|\sigma\|_{C([0,t];H^\alpha)} \leq C \sup_{0\leq z\leq t}\int_0^z \tau^{\gamma/2-1}\|\sigma\|_{C([0,z-\tau];H^\alpha)}d\tau.$$

This implies that $\|\sigma\|_{C([0,t];H^\alpha)} = 0$ for $t \in [0,T]$. The uniqueness is then shown.

(ii).

Assume $\rho(x,t)$ is the mild solution, which satisfies

$$\rho(x,t) = E_\gamma(-t^\gamma A)\rho_0 + \gamma \int_0^t (t-s)^{\gamma/2-1}((t-s)^{\gamma/2}E'_\gamma(-(t-s)^\gamma A)(-\nabla \cdot (\rho a)|_s))\, ds.$$

For the integral, we have done change of variables $s = t - \tau$. Note that $E_\gamma(-t^\gamma A)\rho_0 \in C^\gamma([0,T]; H^\alpha(\mathbb{R}^2)) \cap C^\infty((0,\infty); H^\alpha(\mathbb{R}^2))$ by [25, Equation (8.13)]. Since

$$(t-s)^{\gamma/2}E'_\gamma(-(t-s)^\gamma A)\nabla : H^\alpha(\mathbb{R}^2) \to H^\alpha(\mathbb{R}^2)$$

is a bounded operator and $\rho \in C([0,\infty), H^\alpha(\mathbb{R}^2))$, we apply Proposition 3.3 repeatedly, and find that for any $T > 0$
$$\rho \in C^{0,\gamma}([0,T]; H^\alpha(\mathbb{R}^2)) \cap C^\infty((0,\infty); H^\alpha(\mathbb{R}^2)).$$

Since $E_\gamma(-t^\gamma A)\varphi$ solves the fractional diffusion equation, we have
$$E_\gamma(-t^\gamma A)\varphi = \varphi - \frac{1}{\Gamma(\gamma)} \int_0^t (t-s)^{\gamma-1} A E_\gamma(-s^\gamma A)\varphi \, ds.$$

Secondly, taking the derivative on $t$, we find the operator identity,
$$-\gamma t^{\gamma-1} E'(-t^\gamma A) = -\frac{1}{\Gamma(\gamma)} t^{\gamma-1} I + \frac{\gamma}{\Gamma(\gamma)} \int_0^t (t-s)^{\gamma-1} A s^{\gamma-1} E'(-s^\gamma A) \, ds.$$

Using these two identities and the fact $A\rho \in C([0,\infty); H^{\alpha-2}(\mathbb{R}^2))$, we find that the mild solution satisfies in $C^\gamma([0,T]; H^{\alpha-2}(\mathbb{R}^2)) \cap C^\infty((0,\infty); H^{\alpha-2}(\mathbb{R}^2))$ that

(A.2) $$\rho(x,t) = \rho_0 + \frac{1}{\Gamma(\gamma)} \int_0^t (t-s)^{\gamma-1} (-A\rho(s) - \nabla \cdot (\rho a)(s)) \, ds.$$

Using these time regularity and (A.2), we find
$$D_c^\gamma \rho = \frac{1}{\Gamma(1-\gamma)} \int_0^t \frac{\dot\rho(s)}{(t-s)^\gamma} ds = -\nabla \cdot (\rho a) + \Delta \rho$$
holds in $C([0,T]; H^{\alpha-2}(\mathbb{R}^2))$.

(iii).
For the positivity, it is a little tricky. We first introduce the following notations

(A.3) $$v^+ = \max(v,0), \quad v^- = -\min(v,0).$$

Then, if $v \in H^1(\mathbb{R}^2)$, then $v^\pm \in H^1(\mathbb{R}^2)$ and $\|v^\pm\|_{H^1} \leq \|v\|_{H^1}$.

The idea is then to consider a modified equation
$$D_c^\gamma v = -\nabla \cdot (v^+ a(x,t)) + \Delta v.$$

By the same techniques as in the proof for (i), this equation has a unique global mild solution in $C([0,T]; H^1(\mathbb{R}^2))$ and $v \in C^0([0,T]; H^{-1}(\mathbb{R}^2)) \cap C^1((0,T); H^{-1}(\mathbb{R}^2))$ so that in $C([0,T]; H^{-1}(\mathbb{R}^2))$, we have
$$\frac{1}{\Gamma(1-\gamma)} \int_0^t \frac{\dot v(s)}{(t-s)^\gamma} ds = -\nabla \cdot (v^+ a) + \Delta v.$$

By Proposition 2.18, we have in $C([0,T]; H^{-1}(\mathbb{R}^2))$,
$$\frac{1}{\Gamma(1-\gamma)} \left( \frac{v(t) - v(0)}{t^\gamma} + \gamma \int_0^t \frac{v(t) - v(s)}{(t-s)^{\gamma+1}} ds \right) = -\nabla \cdot (v^+ a(x,t)) + \Delta v.$$

Since $H^{-1}(\mathbb{R}^2)$ is the dual space of $H^1(\mathbb{R}^2)$, we can multiply $v^- = -\min(v,0) \geq 0$ which is in $H^1(\mathbb{R}^2)$ and integrate,
$$\Gamma(1-\gamma)\|\nabla v^-\|_2^2 = \Big(\frac{-\|v^-\|_2^2 - \int \rho_0 v^- \, dx}{t^\gamma}$$
$$+ \gamma \int_0^t \frac{-\|v^+(s)v^-(t)\|_1}{(t-s)^{\gamma+1}} ds - \gamma \int_0^t \frac{\int (v^-(t) - v^-(s))v^-(t) \, dx}{(t-s)^{\gamma+1}} ds\Big)$$
$$\leq \left( \frac{-\|v^-\|_2^2}{2t^\gamma} - \gamma \int_0^t \frac{\int (v^-(t) - v^-(s))v^-(t) dx}{(t-s)^{\gamma+1}} ds \right).$$

Further, noting that $-(v^-(t) - v^-(s))v^-(t) \leq -\frac{1}{2}((v^-(t))^2 - (v^-(s))^2)$, we have

$$\|\nabla v^-\|_2^2 \leq -\frac{1}{2}D_c^\gamma \|v^-\|_2^2,$$

or

$$\frac{1}{2}D_c^\gamma \|v^-\|_2^2 \leq -\|\nabla v^-\|_2^2.$$

By Lemma 2.3, we find that $v^- = 0$. This means that $v$ also solves the original equation and thus $v = \rho$ a.e. by the uniqueness of mild solutions. Hence, $\rho \geq 0$.

For the mass conservation, we consider the approximation sequence (A.1) again. Note that $\rho \in C([0,T]; L^2)$. Let $Y$ be defined in [29, eq. 2.20]. We have

$$\left\| \int_0^t \tau^{\gamma-1} E_\gamma'(-\tau^\gamma A) \nabla \cdot (\rho^{n-1} a)|_{t-\tau} \, d\tau \right\|_1$$
$$\leq \int_0^t \left\| \int_{\mathbb{R}^2} \nabla Y(x-y,\tau) \cdot (\rho^{n-1} a)(y, t-\tau) \, dy \right\|_1 d\tau.$$

Using the estimates for $\nabla Y$ in [29, Lemma 4.7], we have $\|\nabla Y(\cdot, t)\|_1 \leq C t^{\frac{\gamma}{2}-1}$. We therefore find that

$$\|\rho^n\|_1(t) \leq \|\rho_0\|_1 + \frac{\alpha}{\Gamma(\gamma/2)} \int_0^t \tau^{\frac{\gamma}{2}-1} \|\rho^{n-1}(\cdot, t-\tau)\|_1 \, d\tau.$$

for some constant $\alpha$. Consequently, $\|\rho^n\|_1$ is controlled by the solution to $D_c^{\gamma/2} u = \alpha u$ with initial value $\|\rho_0\|_1$. Hence, we conclude that $\|\rho\|_1(t) \leq u(t)$ as well. Consider (5.25). It is easy to see

$$\int_{\mathbb{R}^2} E_\gamma(-t^\gamma A) \rho_0 \, dx = \int_{\mathbb{R}^2} \rho_0 \, dx.$$

The integrand in the second term is integrable since $\|\rho\|_1(t) \leq u(t)$. Then, by Fubini, we can integrate in $x$ first and thus have

$$\int_{\mathbb{R}^2} \int_0^t \tau^{\gamma-1} E_\gamma'(-\tau^\gamma A) \nabla \cdot (\rho a)|_{t-\tau} \, d\tau dx = \int_0^t \int_{\mathbb{R}^2} \tau^{\gamma-1} E_\gamma'(-\tau^\gamma A) \nabla \cdot (\rho a)|_{t-\tau} \, dx d\tau.$$

For $\tau > 0$, $\nabla Y$ is integrable, and its integral must be zero. Hence, for almost every $\tau \in (0, t)$,

$$\int_{\mathbb{R}^2} \nabla E_\gamma'(-\tau^\gamma A) \cdot (\rho a|_{t-\tau}) \, dx = 0.$$

The second term is therefore zero. □